\documentclass{article}

\usepackage[final]{corl_2026} 

\usepackage{amsmath,amsfonts,bm}
\usepackage{mathtools}
\usepackage{physics}

\def\eqref#1{equation~\ref{#1}}

\def\1{\bm{1}}

\DeclareMathAlphabet{\mathsfit}{\encodingdefault}{\sfdefault}{m}{sl}
\SetMathAlphabet{\mathsfit}{bold}{\encodingdefault}{\sfdefault}{bx}{n}

\newcommand{\E}{\mathbb{E}}

\newcommand{\R}{\mathbb{R}}

\DeclareMathOperator*{\minimize}{minimize}
\DeclareMathOperator*{\maximize}{maximize}
\DeclareMathOperator*{\st}{subject\ to}

\DeclareMathOperator{\diag}{diag}
\DeclareMathOperator{\blkdiag}{blkdiag}
\DeclareMathOperator{\nnz}{nnz}
\DeclareMathOperator{\sign}{sign}

\DeclareMathOperator{\detach}{detach}

\renewcommand{\norm}[1]{\lVert#1\rVert}

\newcommand{\defeq}{\vcentcolon=}

\newcommand{\hquad}{\hspace{0.5em}}

\usepackage[capitalize]{cleveref}
\usepackage{booktabs}
\usepackage{caption}

\makeatletter
\newcommand{\smallersize}{\@setfontsize\smallersize\@viiipt\@ixpt}
\makeatother

\title{Scalable Deep Unfolding of Conic Optimizers}

\author{
  Alex~Oshin~\quad~Rahul~Vodeb~Ghosh~\quad~Evangelos~A.~Theodorou\\
  Georgia~Institute~of~Technology\\
  \texttt{\{alexoshin,rghosh88,evangelos.theodorou\}@gatech.edu}
}

\begin{document}
\maketitle


\begin{abstract}
  Deep unfolding (DU) accelerates iterative optimizers by training learnable components through unrolled iterations, but extending DU to the large-scale semidefinite programs (SDPs) common in robotics has remained limited. Unrolling a full-update conic solver such as COSMO exposes two obstacles: backpropagating through the per-iteration linear-system solve incurs memory quadratic in the problem size when the coefficient matrix is formed explicitly, and backpropagating through the positive semidefinite (PSD) cone projection is numerically unstable under repeated eigenvalues. We address the first with a matrix-free implicit differentiation rule that reduces memory from $O(n^2)$ to $O(n + \mathrm{nnz}(A))$, and the second with a stable backward rule based on the Dale\v{c}kii--Krein representation of the Fr\'echet derivative. Together, these enable learning lightweight hyperparameter policies and warm starts for COSMO. On standalone SDPs and second-order cone programs, ranging from max-cut and Lov\'asz $\vartheta$ to robust estimation and control, the learned policies outperform tuned COSMO on every problem class, reducing iterations to tolerance by up to 60$\times$, and transfer without retraining to instances up to 10$\times$ larger in the size parameter while retaining 2--5$\times$ speedups. Our matrix-free solver implementation\footnote{Code is available at: \url{https://github.com/alexoshin/deep-conic-optimization}} also batches naturally on the GPU, solving 64 Lov\'asz $\vartheta$ SDPs 14$\times$ faster than sequential COSMO.jl. When embedded in sequential convex programming for nonlinear covariance steering, our learned optimizers outperform the interior-point solver Clarabel in terms of both success rate and solve time.
\end{abstract}



\section{Introduction}

Deep unfolding (DU) is a model-based machine learning technique that enhances iterative optimization algorithms through data-driven supervised learning~\citep{Hershey2014Deep,Monga2021Algorithm}, with state-of-the-art results across a variety of problem classes ranging from compressed sensing~\citep{Chen2021Hyperparameter} and video reconstruction~\citep{DeWeerdt2024Deep} to the distributed~\citep{Saravanos2025Deep} and nonlinear programming~\citep{Oshin2026Deep} problems that are commonly found in robotics applications.
In DU, a specific optimizer is first chosen, e.g., gradient descent or the alternating direction method of multipliers (ADMM)~\citep{Boyd2010Distributed}. Learnable components are then introduced into the optimizer equations, such as the step size of gradient descent or the penalty parameters of ADMM.
Finally, the optimizer is unrolled for some fixed number of iterations and is trained using backpropagation through time to minimize a meta-level objective that encourages the desired final properties of the learned optimizer, such as fast convergence or robustness to disturbances.

A key feature of DU is that the number of unrolled iterations can be fixed a priori based on the available computational budget, and the optimizer is then trained to maximize performance within that budget.
This makes DU particularly relevant for real-time robotics problems in which a decision must be executed within a fixed time, e.g., model predictive control (MPC).
Indeed, this property was key for popularizing the method in the signal and image processing communities~\citep{Gregor2010Learning,Wang2015Deep,Balatsoukas-Stimming2019Deep}, and DU now constitutes a leading approach for sparse recovery and computational imaging~\citep{Liu2019ALISTA,Chen2021Hyperparameter,DeWeerdt2024Deep,Wang2025Proximal}.

Recent work applies DU to quadratic and conic solvers by learning warm starts or hyperparameters.
\citet{Sambharya2023EndToEnd,Sambharya2024LearningWarmStart} learn warm-start networks for OSQP~\citep{Stellato2020OSQP} and SCS~\citep{ODonoghue2016Conic} by unrolling their fixed-point iterations, and \citet{Xiong2025Solving} learn SCS warm starts with an unrolled Douglas--Rachford scheme (DR-GD) that replaces each linear-system solve with a gradient step.
Since large warm-start networks tend to be expensive, \citet{Sambharya2024LearningAlgorithm} instead learn lightweight hyperparameter trajectories, and \citet{Saravanos2025Deep} extend the approach to learn feedback policies for OSQP and distributed OSQP.
\citet{Oshin2026Deep} apply DU to sequential convex programming (SCP), the setting we adopt in this work, and introduce a hyperparameter policy parameterization that we extend to non-separable cones such as the second-order cone (SOC) and the positive semidefinite (PSD) cone.
None of these approaches scale to the large-scale semidefinite programs (SDPs) common in robotics, where unrolling a full-update ADMM solver such as COSMO~\citep{Garstka2021COSMO} surfaces two distinct obstacles --- one in the linear-system update and one in the conic projection. We provide an extended related work section in the appendix.

The first obstacle arises when backpropagating through the linear-system solve.
Each iteration of COSMO solves a sparse linear system, which the conjugate gradient (CG) method handles using only matrix-vector products that batch well on GPUs, and the adjoint system for the gradient shares the \emph{same} coefficient matrix.
However, this matrix contains $A^\top A$, where $A$ is the constraint matrix, which can be dense even when $A$ is sparse, so forming it explicitly, as in the implicit differentiation scheme of \citet{Saravanos2025Deep}, incurs $O(n^2)$ memory.
DR-GD~\citep{Xiong2025Solving} avoids differentiating the solve, but its linearized iterations converge substantially slower than full ADMM.
We instead retain the full ADMM update and differentiate the solve implicitly using only matrix-vector products.

The second obstacle arises in the projection onto the PSD cone, which requires an eigenvalue decomposition to clamp the negative eigenvalues to zero.
While this forward operation is well-posed, backpropagating through the eigendecomposition is undefined when eigenvalues coincide, owing to a term $\lambda_i - \lambda_j$ appearing in the denominator of the gradient.
We observe that repeated eigenvalues occur frequently during training, particularly under the isotropic covariance targets common in our setting, so resolving this is essential for learning to solve SDPs.
Building on the work of \citet{Engin2018DeepKSPD}, we apply the Dale\v{c}kii--Krein representation of the Fr\'echet derivative~\citep{Higham2008Functions} as a numerically stable backward rule for the PSD projection.

\begin{figure}[h]
  \centering
  \includegraphics[width=\linewidth]{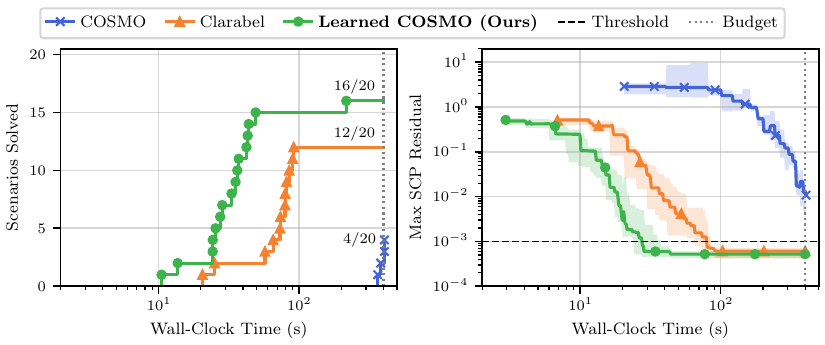}
  \caption{%
    Learned COSMO as the subproblem solver in sequential convex programming (SCP) for nonlinear unicycle covariance steering, over 20 randomized scenarios.
    A run succeeds when the max SCP residual falls below $10^{-3}$ (dashed), within a budget of 20 SCP iterations and up to 20\,s per subproblem solve (dotted, 400\,s).
    Failure modes and closed-loop validation are provided in \cref{appendix:control_performance}.
    Clarabel~\citep{Goulart2024Clarabel} reuses its factorization across SCP iterations.
    \emph{Left:} Scenarios solved vs.\ wall-clock time.
    \emph{Right:} Max SCP residual on the 12 scenarios solved by both learned COSMO and Clarabel (a subset of ours; median and IQR).
    Learned COSMO reaches the threshold in a median of 27.9\,s vs.\ 76.7\,s for Clarabel, and is faster on every one of these scenarios.
  }
  \label{fig:nonlinear_cs_comparison}
\end{figure}

\clearpage
In summary, this paper makes the following contributions:
\begin{itemize}
  \item A matrix-free implicit-differentiation rule for the linear-system solve that avoids instantiating the dense $A^\top A$ coefficient matrix, reducing memory from $O(n^2)$ to $O(n + \mathrm{nnz}(A))$ and enabling backpropagation at scales where direct factorization runs out of memory.
  \item A numerically stable backward rule for the PSD-cone projection based on the Dale\v{c}kii--Krein Fr\'echet derivative, which remains well-defined even under the repeated eigenvalues that arise frequently during the training process.
  \item A demonstration that, together, these enable unrolling a full-update conic solver (COSMO) on large-scale SDPs and SOCPs, with the learned solver requiring up to 60$\times$ fewer iterations than tuned COSMO (\cref{tab:cosmo_policies}).
  The learned policies transfer without retraining to instances up to 10$\times$ larger in the size parameter (\cref{tab:ood_generalization}), and batched GPU solves are up to 14$\times$ faster than COSMO.jl (\cref{tab:batched_solves}).
  \item An application to chance-constrained covariance steering (CS) with nonlinear unicycle dynamics via SCP, where the learned solver succeeds on 16/20 scenarios, compared with 4/20 for COSMO and 12/20 for Clarabel (\cref{fig:nonlinear_cs_comparison}).
\end{itemize}


\section{Methodology}
\label{sec:method}

We consider the deep unfolding of a full-update ADMM-based conic solver. We first define the problem and notation, then identify the two per-iteration steps whose gradients are problematic at scale --- the linear-system solve and the conic projection --- which motivate the two contributions developed in \cref{sec:matrix_free,sec:psd_backward}. We focus on COSMO~\citep{Garstka2021COSMO} as the concrete solver, but emphasize that these gradient rules apply to any ADMM-based conic solver sharing these steps (e.g., SCS~\citep{ODonoghue2016Conic}).

\subsection{Deep Unfolding Conic Optimization}
\label{sec:unfolding-setup}

We consider convex conic programs in the standard form
\begin{equation} \label{eq:conic-standard}
  \minimize_x \hquad \tfrac{1}{2} x^\top P x + q^\top x, \quad \st \hquad A x + s = b, \hquad s \in \mathcal{K},
\end{equation}
where $P \in \mathbb{S}^n_+$, $q \in \R^n$, $A \in \R^{m \times n}$, $b \in \R^m$, and $\mathcal{K}$ is a Cartesian product of convex cones (e.g., nonnegative orthant, SOC, PSD cone).
Given an instance $d = (P, q, A, b, \mathcal{K})$, an ADMM solver produces a sequence of iterates $z^{0}, z^{1}, \ldots$ via a fixed-point map $z^{k+1} = \mathcal{T}_\theta(z^{k}; d)$ parameterized by algorithm hyperparameters $\theta$ (e.g., the penalty parameters $\sigma, \rho$).

In deep unfolding, we make $\theta$ (and, optionally, the initial iterate $z^{0}$) the output of a learnable map $\theta = \pi_\phi(\cdot)$ with parameters $\phi$, unroll the solver for a fixed budget of $K$ iterations, and train $\phi$ by backpropagation through time to minimize a meta-objective $\mathcal{L}$ that rewards rapid convergence:
\begin{equation} \label{eq:du-objective}
  \minimize_{\phi} \hquad \mathbb{E}_{d \sim \mathcal{D}}\big[\, \mathcal{L}(z^{0}, \ldots, z^{K}) \,\big].
\end{equation}

\subsection{Conic Operator Splitting Method (COSMO)}
\label{sec:cosmo-background}

COSMO is a first-order method for large-scale conic programs, such as SDPs~\citep{Garstka2021COSMO}.
Each ADMM iteration requires the solution of an equality-constrained quadratic program (QP), followed by a projection onto the convex cone $\mathcal{K}$, and ends with a dual variable update:
\begin{subequations}
\begin{equation}
  \begin{bmatrix} P + \sigma I & A^\top \\ A & -D(\rho^k)^{-1} \end{bmatrix} \begin{bmatrix} \tilde{x}^{k + 1} \\ \tilde{\nu}^{k + 1} \end{bmatrix} = \begin{bmatrix} \sigma x^k - q \\ b - s^k + y^k \oslash \rho^k \end{bmatrix}, \label{eq:cosmo_linear_system}
\end{equation}
\begin{align}
  \tilde{s}^{k + 1} &= s^k - (y^k + \tilde{\nu}^{k + 1}) \oslash \rho^k, \\
  x^{k + 1} &= \alpha^k \tilde{x}^{k + 1} + (1 - \alpha^k) x^k, \\
  s^{k + 1} &= \Pi_\mathcal{K}(\alpha^k \tilde{s}^{k + 1} + (1 - \alpha^k) s^k + y^k \oslash \rho^k), \label{eq:cosmo_projection} \\
  y^{k + 1} &= y^k + \rho^k \odot (\alpha^k \tilde{s}^{k + 1} + (1 - \alpha^k) s^k - s^{k + 1}),
\end{align}
\end{subequations}
where $\sigma > 0$, $\rho^k \in \R_{++}^m$, and $\alpha^k \in (0, 2)$ are hyperparameters of the algorithm, and $\odot/\oslash$ denote elementwise product and division, respectively.
Our learning framework parameterizes $\rho^k, \alpha^k = \pi_\phi(z^k)$ as a small feedforward network and trains it by minimizing~\cref{eq:du-objective}.
We fix $\sigma = 10^{-6}$ following~\citet{Garstka2021COSMO}.
Notably, the algorithm steps consist entirely of scalar-vector and vector-vector operations that can be differentiated easily through automatic differentiation (``autograd''), except for the linear-system solve~\cref{eq:cosmo_linear_system} and the projection~\cref{eq:cosmo_projection}.
These two steps are the main bottlenecks preventing deep unfolding of large-scale conic optimizers, and addressing them is the focus of the remainder of this section.

\subsection{Matrix-Free Implicit Differentiation}
\label{sec:matrix_free}

COSMO offers two modes for solving the linear system~\cref{eq:cosmo_linear_system}.
The direct mode solves the saddle-point system as written through a sparse $LDL^\top$ factorization, exactly as in OSQP~\citep{Stellato2020OSQP}.
The indirect mode eliminates $\tilde{\nu}^{k + 1}$ and solves the reduced normal-equations form
\begin{equation} \label{eq:cosmo_linear_system_reduced}
  \big( P + \sigma I + A^\top D(\rho^k) A \big) \tilde{x}^{k + 1} = \sigma x^k - q + A^\top (y^k + \rho^k \odot (b - s^k))
\end{equation}
using the conjugate gradient (CG) method, which requires only matrix-vector products (``matvecs'') with the coefficient matrix $M(\rho^k) \defeq P + \sigma I + A^\top D(\rho^k) A$.
The indirect mode is preferred when $n + m$ is large and for GPU implementations~\citep{Schubiger2020GPU}, precisely because it is matrix-free in the forward pass and the CG matvecs batch efficiently.
We adopt the indirect mode in this work, as the unrolled, batched, GPU setting of~\cref{eq:du-objective} is exactly the regime for which it is intended.

\Cref{eq:cosmo_linear_system_reduced} defines the solution $\tilde{x}$ as an implicit function of the matrix $M$ and the right-hand side $\xi \defeq \sigma x^k - q + A^\top (y^k + \rho^k \odot (b - s^k))$.
Using autograd to backpropagate through the CG iterates that produced $\tilde{x}$ scales poorly, as it requires storing and differentiating each CG iteration.
We instead differentiate through \cref{eq:cosmo_linear_system_reduced} directly using implicit differentiation~\citep{Amos2019Differentiable}.
Notably, prior work uses a similar approach to differentiate $M(\rho) \tilde{x} = \xi$ with respect to the \emph{entire} matrix $M$~\citep{Saravanos2025Deep}, and computes $dM/d\rho$ using autograd. This requires $O(n^2)$ memory, despite the fact that only gradients with respect to $\rho$ (a vector of size $m$) are required through the operator.
Therefore, we design a matrix-free implicit differentiation rule that avoids forming these intermediate matrix derivatives.
Note that $\rho$ affects $\tilde{x}$ through both the matrix $M(\rho)$ and the right-hand side $\xi(\rho)$.
Our custom rule handles only the dependence through $M(\rho)$, while the dependence through $\xi(\rho)$ is handled via autograd outside of the custom layer.
By the chain rule, the total derivative $d\tilde{x} / d\rho$ is the sum of these two contributions, so the resulting gradient is exact.

To start, let $\mathcal{L}(\tilde{x})$ be some downstream loss and consider the derivative with respect to a single $\rho_i$:
\begin{equation} \label{eq:grad_rho}
  \pdv{\mathcal{L}}{\rho_i} = \nabla_{\tilde{x}} \mathcal{L}^\top \pdv{\tilde{x}}{\rho_i} = -\nabla_{\tilde{x}} \mathcal{L}^\top M^{-1} \pdv{M}{\rho_i} \tilde{x} = \tilde{\lambda}^\top \pdv{M}{\rho_i} \tilde{x},
\end{equation}
where $\tilde{\lambda}$ is defined as the solution to the adjoint system $M \tilde{\lambda} = -\nabla_{\tilde{x}} \mathcal{L}$, which can be solved for efficiently using the same CG method as the forward pass, but with a different right-hand side.
The same adjoint also gives the gradient with respect to the right-hand side, $\nabla_{\xi} \mathcal{L} = M^{-1} \nabla_{\tilde{x}} \mathcal{L} = -\tilde{\lambda}$, which the custom layer returns so that autograd can propagate the dependence through $\xi(\rho)$ at no additional cost.
Since $M$ only depends on $\rho_i$ through the term $A^\top D(\rho) A = \sum_j \rho_j a_j a_j^\top$, we have $\partial M / \partial \rho_i = a_i a_i^\top$.
Substituting into \cref{eq:grad_rho} and stacking yields:
\begin{equation}
  \pdv{\mathcal{L}}{\rho_i} = (A \tilde{\lambda})_i (A \tilde{x})_i \implies \nabla_\rho \mathcal{L} = (A \tilde{\lambda}) \odot (A \tilde{x}).
\end{equation}

This shows that the gradient can be computed using matvecs alone.
The dense $A^\top A$ term is never formed explicitly in either the forward or the backward pass, so the overall complexity is $O(n + \mathrm{nnz}(A))$ rather than $O(n^2)$. We implement this approach as a custom backward rule for the linear-solve layer and validate its scaling in \cref{fig:linear_system_comparison}. Full experimental details are described in~\cref{appendix:linear_system_experiment}.

\begin{figure}[h]
  \centering
  \includegraphics[width=\linewidth]{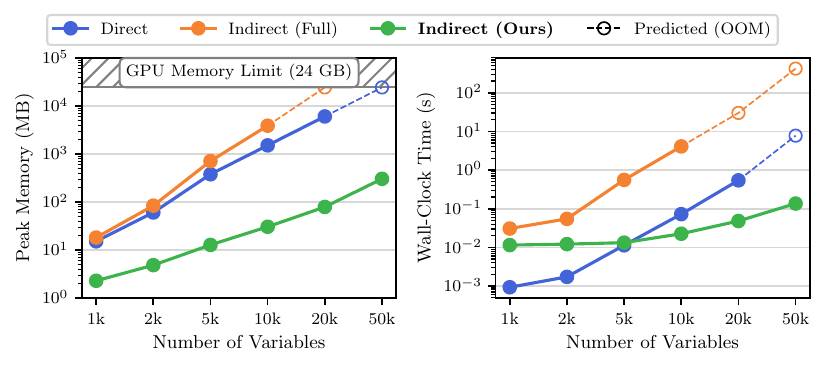}
  \caption{%
    Backward pass cost of the COSMO linear-system solve as a function of problem size.
    We differentiate the loss $\tfrac{1}{2}\norm{x - x^\star}^2$ using three methods:
    autograd through the direct mode~\cref{eq:cosmo_linear_system}, autograd through the indirect mode~\cref{eq:cosmo_linear_system_reduced} using the approach from \citet{Saravanos2025Deep}, and our proposed matrix-free rule.
    \emph{Left:} Peak GPU memory. \emph{Right:} Backward wall-clock time, both versus the number of variables $n$ (with $m = n/2$ constraints), on a single RTX~4090 (24\,GB) in \texttt{float64}.
    Open markers on dashed lines extrapolate methods that exceed the 24\,GB budget.
    Both baselines grow quadratically and exhaust memory by $n \approx 20$--$50$k, whereas our method scales linearly, staying under 300\,MB and remaining the fastest at the largest sizes.
  }
  \label{fig:linear_system_comparison}
\end{figure}

\subsection{Stable PSD Projection}
\label{sec:psd_backward}

\begin{figure}[h]
  \centering
  \includegraphics[width=\linewidth]{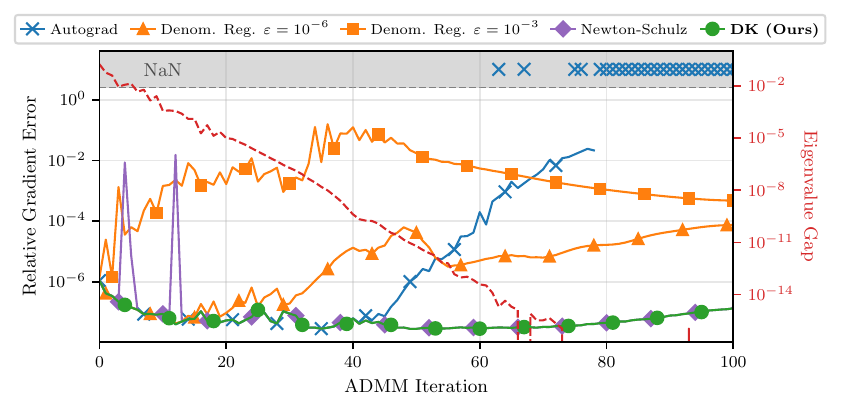}
  \caption{%
    PSD-projection gradient along a real ADMM trajectory. We unroll COSMO on a chance-constrained covariance-steering SDP (double-integrator, $N=10$) from a random initial iterate.
    As the iterates approach the PSD cone, the minimum eigenvalue gap $\min_{i \neq j} \lvert \lambda_i - \lambda_j \rvert$ of $S_k$ (dashed, right axis) falls by over ten orders of magnitude, causing each baseline to eventually fail.
    Autograd's error grows as the gap shrinks until it returns NaNs (top band).
    Denominator regularization is biased whenever the gap drops below $\varepsilon$, so a smaller $\varepsilon$ only delays the error.
    Newton--Schulz produces intermittent error spikes.
    DK stays at the error floor for the whole trajectory and requires no tuning.
  }
  \label{fig:psd_projection_comparison_1}
\end{figure}

The second key operation that requires custom handling arises during the backward pass of the projection onto $\mathcal{K}$ in \cref{eq:cosmo_projection}, specifically when $\mathcal{K}$ involves the PSD cone.
For the zero, nonnegative-orthant, and second-order cones, the projection is smooth (or piecewise-smooth) and its backward pass is stable.
On the other hand, projecting a symmetric matrix $S$ onto the PSD cone requires an eigendecomposition $S = U \Lambda U^\top$ followed by clamping the negative eigenvalues:
\begin{equation} \label{eq:psd-proj}
  \Pi_{\mathbb{S}_+}(S) = U \max(\Lambda, 0) U^\top .
\end{equation}
While this forward projection is well-defined, the backward pass differentiates the full eigendecomposition, which is unstable when eigenvalues become nearly degenerate.
Given cotangents $\bar{U}$ and $\bar{\Lambda} = \diag(\bar{\lambda}_1, \ldots, \bar{\lambda}_n)$, the cotangent $\bar{S}$ is given by:
\begin{equation} \label{eq:psd-cotangent}
  \bar{S} = U \left( \bar{\Lambda} + F \odot (U^\top \bar{U}) \right) U^\top,
\end{equation}
where $F$ is zero on the diagonal and $F_{ij} = 1 / (\lambda_j - \lambda_i)$ for $i \neq j$.
These off-diagonal terms blow up as eigenvalues coalesce and are undefined when they coincide.
We emphasize that this phenomenon is not pathological and occurs routinely in our setting.
For example, in covariance steering problems, we often want to enforce isotropic covariance constraints of the form $S = \sigma^2 I$.
As the constraint gets closer to being satisfied, all the eigenvalues approach $\sigma^2$ asymptotically, causing numerical instability and poor gradients.
This is demonstrated on a real covariance steering problem in \cref{fig:psd_projection_comparison_1}.

The projection~\cref{eq:psd-proj} is a primary matrix function of the form $T = \Pi_{\mathbb{S}_+}(S) = U f(\Lambda)\, U^\top$, with $f(x) = \max(x, 0)$.
Functions of this form have a Fr\'echet derivative given by the Dale\v{c}kii--Krein (DK) theorem~\citep[Theorem 1]{Engin2018DeepKSPD}.
Given the symmetric cotangent $\bar{T}$ of the output of the projection, the cotangent of the input matrix $\bar{S}$ is given by:
\begin{equation} \label{eq:psd-cotangent-dk}
  \bar{S} = U(L \odot (U^\top \bar{T} U)) U^\top ,
\end{equation}
where $L$ is the L\"owner matrix of first divided differences of $f$ evaluated at the eigenvalues:
\begin{equation} \label{eq:lowner}
  L_{ij} = \begin{cases}
    \frac{f(\lambda_j) - f(\lambda_i)}{\lambda_j - \lambda_i}, & \lambda_i \neq \lambda_j, \\
    f'(\lambda_i), & \lambda_i = \lambda_j.
  \end{cases}
\end{equation}

Contrasting this rule with the one in~\cref{eq:psd-cotangent}, the factor $1 / (\lambda_j - \lambda_i)$ is canceled by a matching numerator $f(\lambda_j) - f(\lambda_i)$ that vanishes at the same rate, and the singularity at $\lambda_i = \lambda_j$ is removable with limit $f'(\lambda_i)$.
For $f(x) = \max(x, 0)$, we have that $f'(\lambda_i) \in [0, 1]$, so the derivative is always finite.
Autograd, by contrast, differentiates the eigenvectors separately, so the reciprocal term is not multiplied by a term that vanishes as eigenvalues coalesce, and the derivative diverges.
For repeated zero eigenvalues, we adopt the \texttt{torch.relu} subgradient convention $f'(0) = 0$.

Crucially, this DK rule can be implemented for free in most modern deep learning libraries.
Since the forward pass already requires computing the eigendecomposition $S = U \Lambda U^\top$, the backward pass can simply reuse the decomposition to compute $\bar{S}$ at the cost of one Hadamard product and four matrix multiplications, comparable to the na\"ive autograd backward pass it replaces, which must additionally backpropagate through $U f(\Lambda) U^\top$ to obtain the cotangents in~\cref{eq:psd-cotangent}.

\paragraph{Alternative Approaches.}

\begin{figure}[h]
  \centering
  \includegraphics[width=\linewidth]{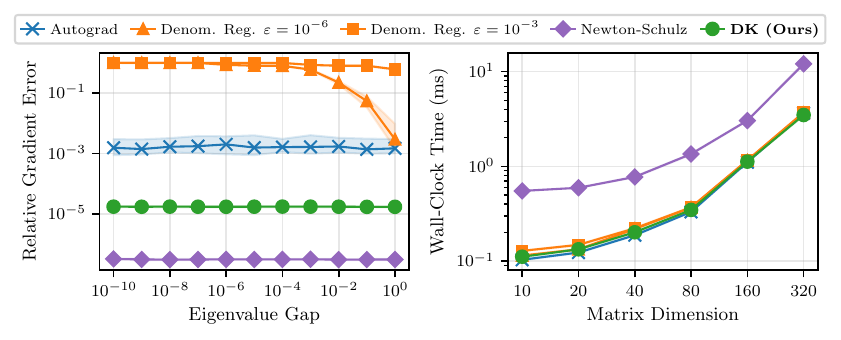}
  \caption{%
    DK is free in the regime where autograd survives. \emph{Left:} Relative error of each backward rule against a finite-difference reference on matrices of size $n = 40$ with clustered eigenvalues. Autograd remains finite, but is less accurate than DK and Newton--Schulz. Denominator regularization carries a large bias, even when the gap $\delta \gg \varepsilon$, because the gaps within each cluster ($\sim 10^{-3}\delta$) remain small. \emph{Right:} Wall-clock time of the backward rule vs.\ matrix dimension. Autograd, denominator regularization, and DK have comparable costs, whereas Newton--Schulz's unrolled iterations make it nearly an order of magnitude more expensive.
  }
  \label{fig:psd_projection_comparison_2}
\end{figure}

An initial idea to handle this repeated-eigenvalue singularity would be to perturb the $\lambda_j - \lambda_i$ term appearing in the denominator of $F_{ij}$ with some small regularization that numerically avoids the division by zero.
We implemented the ``denominator regularization'' rule
\begin{equation} \label{eq:denom_reg}
  F_{ij} = \frac{\lambda_j - \lambda_i}{(\lambda_j - \lambda_i)^2 + \varepsilon^2},
\end{equation}
which returns zero at an exact tie and reduces to $1 / (\lambda_j - \lambda_i)$ for gaps much larger than $\varepsilon$.
It also recovers DK as $\varepsilon \to 0$.
While this rule is numerically stable, it results in a large bias when $\varepsilon$ is greater than the true gap (\cref{fig:psd_projection_comparison_2}), and this affects the quality of the downstream gradients (\cref{fig:psd_projection_comparison_1}).

Another potential approach avoids the eigendecomposition gradient entirely by re-expressing the projection using the matrix sign function.
Using $\Pi_{\mathbb{S}_+}(S) = \tfrac{1}{2}(S + |S|)$ with $|S| = \sign(S) S$, one computes $\sign(S)$ by the scaled Newton--Schulz iteration
\begin{equation}
  S_{k+1} = \tfrac{1}{2} S_k \big(3I - S_k^2\big),
  \qquad S_0 = S / \lVert S \rVert_F ,
  \label{eq:newton-schulz}
\end{equation}
which converges quadratically to $\sign(S)$ once $\lVert S_k \rVert_2 < 1$.
The backward pass then involves autograd through the $K$ unrolled matrix products, which contain no $1/(\lambda_i - \lambda_j)$ term and so cannot diverge.
This approach is numerically stable, but the main drawback is the computational cost of unrolling $K$ Newton--Schulz iterations during the backward pass (see \cref{fig:psd_projection_comparison_2} for comparison).

\section{Experiments}
\label{sec:results}

We evaluate our DU methodology on nonlinear covariance steering problems, showing how learned conic optimization can be used effectively as a subroutine in SCP methods, and on a suite of standalone SDPs and SOCPs that quantify how well the method generalizes.
We use the policy parameterization from~\citet{Oshin2026Deep}, which is a small 32$\times$32 MLP that maps the current iterate and residuals to the hyperparameters.
For non-separable cones (SOC, PSD), we apply a rectification step (averaging) so that a single penalty parameter is predicted for each cone and the projection remains well-defined.
We then unroll COSMO for 100--400 iterations and train on a dataset of 100 problems per class using DU to minimize a self-supervised loss based on the log of the primal and dual residuals:
\begin{equation}
  \mathcal{L}(z^{0}, \ldots, z^{K}) = \sum_{k = 1}^{K} \log\left( \frac{\norm{Px^k + q + A^\top y^k}^2}{\detach(\norm{Px^0 + q + A^\top y^0}^2)} + \frac{\norm{Ax^k + s^k - b}^2}{\detach(\norm{Ax^0 + s^0 - b}^2)} \right),
\end{equation}
which encourages maximizing the log-convergence rate of the dominant residual~\citep{Nocedal2006Numerical}.
The detached denominators prevent the model from learning a degenerate solution.
To ensure a fair comparison, the non-learned COSMO baselines (tuned constant and tuned adaptive residual-balancing $\rho$) use our own PyTorch reimplementation, run on the same GPU as the learned variant, so that the reported statistics reflect the learned policy itself rather than differences in implementation or hardware.
All experiments were run on a machine with an i9-13900K CPU and an RTX 4090 GPU (warmed up, synchronized timing).
Further training, policy parameterization, and dataset generation details, along with ablations isolating the individual learned components, are provided in the appendix.

\subsection{Nonlinear Covariance Steering}
\label{sec:cs-steering}

For the CS problems, we use an affine feedback parameterization~\citep{Goulart2006Optimization}, which for linear covariance steering results in a convex SDP~\citep{Balci2021Covariance}.
We embed this in a successive convexification scheme inspired by iterative CS and SCP~\citep{Ridderhof2019Nonlinear,Ridderhof2022ChanceConstrained}, linearizing the nonlinear dynamics at each outer iteration to yield a convex SDP that is solved by either COSMO or our learned variant.
We test our approach using both double-integrator and unicycle dynamics.
Since the double integrator is linear, we present the statistics alongside the standalone problems in~\cref{tab:cosmo_policies}.
The learned policy accelerates each inner convex solve, reaching tolerance in a median of 63 iterations vs.\ 3954 for the best-tuned constant COSMO (\cref{tab:cosmo_policies}).
Within SCP, this lets learned COSMO solve 16/20 scenarios vs.\ 4/20 for COSMO, and it is faster than Clarabel on every one of the 12 scenarios that both solve (\cref{fig:nonlinear_cs_comparison}).

\subsection{Benchmark Conic Optimization Problems}
\label{sec:benchmarks}

We also train and benchmark the learned policies on a range of standalone conic programs.
We include max-cut and Lov\'asz $\vartheta$ SDPs on Erd\H{o}s--R\'enyi graphs, which are standard SDP benchmarks.
To cover problems that arise in control and robotics, we include Lyapunov linear matrix inequality (LMI) SDPs for stability certification, and robust Kalman filtering and robust MPC SOCPs for estimation and control under uncertainty.
\Cref{appendix:optimization_problems} gives the full mathematical formulation of each optimization problem.
We emphasize that our scalability contributions are \emph{necessary} for training at these sizes.
For example, on the CS unicycle problems, both the direct and indirect (full) methods from \cref{sec:matrix_free} run out of memory on a 24\,GB GPU, while our matrix-free rule peaks at 1.3\,GB and takes 29\,s per epoch.
\Cref{tab:cosmo_policies} reports iterations to tolerance against COSMO using both tuned constant hyperparameters and a tuned adaptive $\rho$ policy.
Each arm is run until the ADMM residuals fall below an absolute tolerance of $10^{-3}$ or until a max budget of 10000 iterations / 60\,s is reached.
Neither constant nor adaptive $\rho$ generalizes well across all classes, occasionally producing failed or extremely slow solves.
The learned policies increase the number of successful solves and reduce median iterations to tolerance on every class.
They also narrow the spread of iteration counts on every class except robust KF, making solve times far more predictable, which is important for real-time settings with a fixed computational budget such as MPC.
When compared against COSMO.jl~\citep{Garstka2021COSMO} and Clarabel~\citep{Goulart2024Clarabel}, our learned optimizers are fastest on the largest classes (Lov\'asz $\vartheta$ and both CS problems), while the CPU solvers excel on the small SOCPs (\cref{appendix:additional_benchmarks}).

\begin{table}[h]
  \centering
  \caption{%
    Comparison of COSMO hyperparameter policies on benchmark conic optimization problems. \emph{Constant} uses the best-tuned fixed values of $\rho$ and $\alpha$; \emph{Adaptive} uses the best-tuned adaptive $\rho$ policy, based on OSQP-style residual balancing; and \emph{Learned} uses our proposed MLP policies, trained via deep unfolding. We report the median [IQR] of the number of iterations needed to reach an absolute tolerance of $\varepsilon = 10^{-3}$ on 100 held-out test instances, along with the number of successful solves. The learned policies outperform both baselines on every problem class.
  }
  \label{tab:cosmo_policies}
  \smallersize
  \begin{tabular}{l ccc ccc}
    \toprule
    & \multicolumn{3}{c}{Successes (out of 100)} & \multicolumn{3}{c}{Iterations} \\
    \cmidrule(lr){2-4} \cmidrule(lr){5-7}
    Problem Class & Constant & Adaptive & Learned & Constant & Adaptive & Learned \\
    \midrule
    \hquad Lyapunov LMI & \textbf{100} & 47 & \textbf{100} & 590 [397, 998] & 788 [384, 1904] & \textbf{235 [182, 342]} \\
    \hquad Max Cut & \textbf{100} & \textbf{100} & \textbf{100} & 212 [200, 224] & 550 [404, 1034] & \textbf{69 [65, 79]} \\
    \hquad Lov\'asz $\vartheta$ & 0 & 77 & \textbf{100} & -- & 1449 [134, 7848] & \textbf{74 [69, 81]} \\
    \hquad Robust KF & \textbf{100} & \textbf{100} & \textbf{100} & 314 [312, 317] & 314 [312, 317] & \textbf{26 [24, 32]} \\
    \hquad Robust MPC & 99 & \textbf{100} & \textbf{100} & 521 [380, 708] & 81 [59, 102] & \textbf{16 [14, 19]} \\
    \hquad Double Integrator & \textbf{100} & 49 & \textbf{100} & 59 [37, 87] & 684 [367, 1534] & \textbf{18 [10, 31]} \\
    \hquad Unicycle & 90 & 0 & \textbf{99} & 3954 [2672, 5306] & -- & \textbf{63 [60, 68]} \\
    \bottomrule
  \end{tabular}
\end{table}

\subsection{Batched GPU Throughput}
Next, we study the throughput of our implementation when many instances must be solved at once.
Since our linear-system solve is matrix-free (\cref{sec:matrix_free}), a batch of $B$ instances can be stacked into a single block-diagonal \emph{super-batched} conic program and solved in parallel with one CG solve per iteration.
The only operation that differs between single and batched solves is the convergence check, which requires a per-instance reduction over the residuals; otherwise, the solver requires minimal host-device synchronization.
\Cref{tab:batched_solves} compares the wall-clock time to solve an entire batch of max-cut or Lov\'asz $\vartheta$ instances against COSMO.jl and Clarabel, which solve the instances sequentially on the CPU.
Increasing the batch size from $B = 16$ to $B = 64$ increases the learned solver's time by only $1.5\times$ on max cut and $3.4\times$ on Lov\'asz $\vartheta$, compared with $4\times$ or more for the sequential solvers.
The advantage grows with batch size: at $B = 64$, the learned solver is $8\times$ (max cut) and $14\times$ (Lov\'asz $\vartheta$) faster than COSMO.jl.
Full experimental details are given in~\cref{appendix:batched_experiment}.

\begin{table}[h]
  \centering
  \caption{%
    Our learned solver (batched, GPU) vs.\ open-source COSMO.jl v0.8.9 and Clarabel v0.11.1 (sequential, CPU), across varying batch sizes. Ten batches are sampled per size. Median [IQR] of wall-clock times (in seconds) to solve the entire batch.
  }
  \label{tab:batched_solves}
  \small
  \begin{tabular}{l cc cc}
    \toprule
    & \multicolumn{2}{c}{Max Cut} & \multicolumn{2}{c}{Lov\'asz $\vartheta$} \\
    \cmidrule(lr){2-3} \cmidrule(lr){4-5}
    Method & $B=16$ & $B=64$ & $B=16$ & $B=64$ \\
    \midrule
    Learned & \textbf{0.248 [0.232, 0.264]} & \textbf{0.372 [0.361, 0.397]} & \textbf{1.53 [1.21, 2.12]} & \textbf{5.14 [3.65, 8.61]} \\
    COSMO.jl & 0.721 [0.625, 0.757] & 3.02 [2.74, 3.28] & 9.63 [9.13, 10.3] & 72.3 [39.8, 98.2] \\
    Clarabel & 3.70 [3.60, 3.76] & 14.9 [14.7, 17.4] & 68.6 [68.4, 69.0] & 276 [275, 277] \\
    \bottomrule
  \end{tabular}
\end{table}

\subsection{Out-of-Distribution Generalization}
Finally, we analyze the robustness of our learned policies to instances $5$--$10\times$ larger in the size parameter than those seen during training.
\Cref{tab:ood_generalization} uses 20 newly sampled instances for each size, rather than the 100-instance test set from \cref{tab:cosmo_policies}, so the $1\times$ statistics differ.
The learned policies solve every instance and retain a $1.76$--$12.5\times$ speedup in both iterations and wall-clock time over the best-tuned baseline at every size.
Across all classes, the speedups decrease as the problem size increases, but for different reasons.
For robust KF, the baseline's iteration count is unchanged while the learned policy's grows, reflecting distribution shift.
On robust MPC, the learned iteration count stays constant and the gap narrows only because the baseline converges faster at longer horizons.
On max cut, both the baseline and learned iteration counts grow as the PSD block becomes larger.
We leave training on mixed-size datasets to future work.
Full details are given in~\cref{appendix:ood_experiment}.

\begin{table}[h]
  \centering
  \scriptsize
  \caption{Out-of-distribution generalization of learned policies trained at $1\times$ and evaluated (without retraining) on new instances $5$--$10\times$. Median [IQR] iterations and wall-clock time to $\varepsilon = 10^{-3}$ over 20 held-out instances (no arm failed on any instance). \emph{Baseline} is the best-tuned configuration for non-learned COSMO. Speedup is the median of per-instance baseline\,$/$\,learned ratios.}
  \label{tab:ood_generalization}
  \setlength{\tabcolsep}{4pt}
  \begin{tabular}{l l cc cc cc}
    \toprule
                     &              & \multicolumn{2}{c}{Iterations}          & \multicolumn{2}{c}{Wall-Clock Time (s)}           & \multicolumn{2}{c}{Speedup} \\
                     \cmidrule(lr){3-4}                        \cmidrule(lr){5-6}                                  \cmidrule(lr){7-8}
    Problem Class    & Size         & Baseline               & Learned                        & Baseline                    & Learned                            & Iterations         & Time         \\
    \midrule
    Max Cut          & $n{=}50\hquad$ ($1\times$)  & 216 [194, 232]   & \textbf{86 [78, 90]}    & 0.41 [0.36, 0.44]   & \textbf{0.20 [0.18, 0.21]} & $2.53\times$ & $2.02\times$ \\
                     & $n{=}250$ ($5\times$) & 432 [420, 437] & \textbf{204 [192, 214]} & 2.65 [2.58, 2.68]   & \textbf{1.35 [1.27, 1.42]}    & $2.14\times$ & $1.96\times$ \\
                     & $n{=}500$ ($10\times$) & 593 [581, 627] & \textbf{316 [308, 328]} & 	6.66 [6.52, 7.04] & \textbf{3.85 [3.75, 3.99]} & $1.91\times$ & $1.76\times$ \\
    \midrule
    Robust KF        & $N{=}100$ ($1\times$) & 318 [314, 320]    & \textbf{25 [24, 37]}    & 4.72 [4.59, 4.80]   & \textbf{0.84 [0.81, 0.94]}    & $12.5\times$ & $5.49\times$ \\
                     & $N{=}500$ ($5\times$) & 320 [317, 321]    & \textbf{48 [39, 60]}    & 6.50 [6.37, 6.55]   & \textbf{1.42 [1.33, 1.55]}    & $6.69\times$ & $4.58\times$ \\
                     & $N{=}1000$ ($10\times$) & 320 [319, 322]    & \textbf{66 [62, 80]}    & 7.06 [7.03, 7.12]   & \textbf{1.91 [1.82, 2.17]}    & $4.82\times$ & $3.67\times$ \\
    \midrule
    Robust MPC       & $N{=}50$ ($1\times$)  & 66 [52, 91]       & \textbf{16 [14, 17]}    & 0.690 [0.477, 0.881]   & \textbf{0.217 [0.207, 0.233]} & $4.16\times$ & $3.16\times$ \\
                     & $N{=}250$ ($5\times$) & 46 [39, 51]       & \textbf{14 [13, 16]}    & 0.536 [0.486, 0.642]   & \textbf{0.214 [0.191, 0.230]} & $3.31\times$ & $2.59\times$ \\
                     & $N{=}500$ ($10\times$) & 43 [33, 56]       & \textbf{14 [14, 17]}    & 0.584 [0.387, 0.633]   & \textbf{0.224 [0.208, 0.253]}    & $2.64\times$ & $2.16\times$ \\
    \bottomrule
  \end{tabular}
\end{table}


\section{Limitations}
\label{sec:limitations}

Our matrix-free rule removes the $O(n^2)$ memory bottleneck of the linear-system solve, but the PSD projection still requires a dense $O(p^3)$ eigendecomposition of each semidefinite block of size~$p$.
The DK rule makes the backward pass through this projection numerically stable but does not reduce its asymptotic cost.
Thus, our scalability gains apply to growth in the number of variables and constraints, rather than in the size of individual conic blocks.
The learned policies also target low-to-moderate accuracy ($\varepsilon \approx 10^{-3}$), which is typical of the real-time SCP and MPC settings we consider.
Because first-order methods converge slowly at high accuracy~\citep{He2012Convergence}, we do not claim that learned COSMO outperforms interior-point solvers such as Clarabel at tight tolerances.
In addition, as with any deep-unfolding method, we train a separate policy for each problem family.
This upfront cost must be amortized over many solves, although it was modest in our experiments: roughly 100 training problems per class.
Finally, although we instantiate the approach on COSMO, both gradient rules apply to any ADMM-based conic solver that shares these steps (e.g., SCS); we leave an empirical demonstration as well as hardware testing to future work.


\clearpage
\acknowledgments{The authors would like to thank the anonymous reviewers for their insightful comments and feedback, which helped strengthen both the analysis and the presentation of this work.}


\bibliography{references}  


\clearpage
\appendix


\section{Related Work}
\label{appendix:related_work}


\paragraph{Learning to optimize and deep unfolding.}
Deep unfolding, also known as \emph{algorithm unrolling}, is a model-based learning-to-optimize approach where the iterations of an optimization algorithm are unrolled into a differentiable compute graph and a small set of learnable components is trained through backpropagation~\citep{Hershey2014Deep,Monga2021Algorithm,Amos2023Tutorial}.
The technique originated in the signal and image processing communities, beginning with the learned ISTA (LISTA) approximation of sparse coding~\citep{Gregor2010Learning}, which was then extended to many subsequent unrolled sparse recovery networks~\citep{Wang2015Deep,Liu2019ALISTA,Balatsoukas-Stimming2019Deep}.
More broadly, DU sits at the structured end of the learning-to-optimize spectrum~\citep{Chen2022Learning}.
Rather than replacing the update rule with a black-box neural network~\citep{Andrychowicz2016Learning,Li2017Learning}, DU retains the known optimizer and learns only a few of its components, preserving the underlying interpretability and theoretical properties of the optimizer, such as convergence and robustness.
This mirrors gradient-based meta-learning, where an outer objective is optimized through a short unrolled inner optimization~\citep{Finn2017Modelagnostic}.
DU now constitutes a leading approach for sparse recovery and computational imaging~\citep{Chen2021Hyperparameter,DeWeerdt2024Deep,Wang2025Proximal}, where a key feature is that the unroll length can be fixed a priori based on a given computational budget and the optimizer is trained to maximize performance within that budget---precisely the property that makes DU appealing for the real-time, budget-constrained robotics settings we target.

\paragraph{Learned QP and conic solvers.}
Recent work brings DU to quadratic and conic solvers under two broad settings: learned warm starts and learned hyperparameters.
Warm-start networks map problem data to an initial iterate for solvers such as OSQP~\citep{Stellato2020OSQP} and SCS~\citep{ODonoghue2016Conic}, training them through unrolled fixed-point iterations~\citep{Sambharya2023EndToEnd,Sambharya2024LearningWarmStart}.
\citet{Xiong2025Solving} replace each inner linear-system solve with a linearized Douglas--Rachford step to avoid differentiating through it.
Meanwhile, in the learned-hyperparameter setting, the per-iteration hyperparameters of an optimizer are trained to maximize some downstream objective, such as convergence speed~\citep{Sambharya2024LearningAlgorithm}.
These hyperparameters can be parameterized using neural networks, which allows for a form of closed-loop adaptation that can greatly accelerate the optimization~\citep{Ichnowski2021Accelerating,Saravanos2025Deep}.
Moreover, embedding accelerated conic optimizers within a sequential convex programming approach can lead to orders-of-magnitude improvements over traditional approaches~\citep{Oshin2026Deep}.
None of these prior approaches scale to the large-scale SDPs we consider, due to obstacles related to the backpropagation of the linear-system solve and conic projection.

\paragraph{Differentiable optimization layers.}
Deep unfolding is distinct from differentiable optimization layers such as \texttt{OptNet}~\citep{Amos2017OptNet} and \texttt{cvxpylayers}~\citep{Agrawal2019Differentiable}, which differentiate the solution of an optimization problem with respect to its problem data by implicitly differentiating the optimality conditions~\citep{Amos2019Differentiable}.
These solver layers have also been embedded into meta-learning pipelines, e.g., differentiable closed-form ridge-regression solvers for few-shot learning~\citep{Bertinetto2019Metalearning}.
Because the resulting derivative is independent of the trajectory taken to the solution, hyperparameters that determine the optimization process, such as the penalty $\rho^k$ and relaxation $\alpha^k$ parameters that we learn in this work, have an identically zero gradient through these layers, making DU the more appropriate tool when the goal is to accelerate the algorithm itself.
Implicit differentiation~\citep{Blondel2022Efficient} underlies differentiable control layers such as differentiable MPC~\citep{Amos2018Differentiable} and its robust extension~\citep{Oshin2024Differentiable}, which are usually embedded as components in an end-to-end planning or control architecture~\citep{Karkus2023DiffStack}.
These approaches are complementary to our setting, where we instead differentiate the solver itself in order to learn how to accelerate the optimization.

\paragraph{First-order conic solvers.}
We build on the operator-splitting family of first-order conic solvers~\citep{Boyd2010Distributed}, including SCS~\citep{ODonoghue2016Conic}, OSQP~\citep{Stellato2020OSQP}, and COSMO~\citep{Garstka2021COSMO}.
Our matrix-free linear-system update inherits the indirect mode that makes these methods amenable to GPU acceleration~\citep{Schubiger2020GPU}.
First-order methods scale to problem sizes that are intractable for interior-point methods, but converge slowly at high accuracy~\citep{He2012Convergence}.
We therefore target the low-to-moderate accuracy regime and benchmark against the interior-point solver Clarabel~\citep{Goulart2024Clarabel} for a high-accuracy reference.

\paragraph{Differentiating through eigendecompositions.}
Our stable PSD-projection rule is based on the Dale\v{c}kii--Krein representation of the Fr\'echet derivative of a primary matrix function~\citep{Higham2008Functions,Bhatia2009Positive}.
Backpropagation through an eigendecomposition is known to be numerically unstable under degenerate eigenvalues, owing to the $1/(\lambda_i - \lambda_j)$ terms that appear when differentiating the eigenvectors~\citep{Ionescu2015Matrix}.
This has motivated a range of remedies in the deep learning literature, such as the power-iteration method of \citet{Wang2019Backpropagationfriendly} that avoids backpropagating through the eigendecomposition altogether.
\citet{Engin2018DeepKSPD} apply the Dale\v{c}kii--Krein form to learn SPD representations for image recognition.
We show that the same approach can be used as a stable backward rule for the PSD cone projection of an unrolled conic solver, where repeated eigenvalues arise routinely during training.


\section{Details on Linear-System Solve Backward Pass Experiment (Figure 2)}
\label{appendix:linear_system_experiment}

In order to validate our proposed linear-system solve backward pass rule, we sample random conic problem instances of size $n \in \{1\text{k}, 2\text{k}, 5\text{k}, 10\text{k}, 20\text{k}, 50\text{k}\}$ with $m = n/2$ constraints.
Each instance uses a diagonal $P$ with $P_{ii} \sim \mathcal{U}(0, 10)$ and $A$ with 1\% nonzero elements sampled from $A_{ij} \sim \mathcal{N}(0, 1)$.
We randomly generate $\rho$ and a right-hand side $b$ with $\rho_i \sim \mathcal{U}(0, 1)$ and $b_i \sim \mathcal{N}(0, 1)$.
Finally, we sample the least-squares target $x^\star$ with $x^\star_i \sim \mathcal{N}(0, 1)$.

In \cref{fig:linear_system_comparison}, we report the mean peak GPU memory and mean backward pass wall-clock time of each method over 100 random instances per problem size.
We compare three differentiation strategies:
\begin{enumerate}
  \item \textbf{Direct} forms the saddle-point coefficient matrix~\cref{eq:cosmo_linear_system} as a dense matrix and backpropagates through \texttt{torch.linalg.solve} using autograd. This requires $O((n + m)^2)$ memory.
  \item \textbf{Indirect (Full)} forms the reduced coefficient matrix $M = P + \sigma I + A^\top D(\rho) A$ from \cref{eq:cosmo_linear_system_reduced} as a dense matrix and backpropagates through the linear-system solve using implicit differentiation, following~\citet[Theorem 2]{Saravanos2025Deep}. This approach also requires instantiating $\partial \tilde{x} / \partial M$ as a dense matrix during the backward pass, costing $O(n^2)$ memory.
  \item \textbf{Indirect (Ours)} solves and backpropagates through~\cref{eq:cosmo_linear_system_reduced} by treating the coefficient matrix as a linear operator, where $P$ and $A$ are sparse. The backward pass uses the analytic rule \cref{eq:grad_rho}, requiring only matvecs and reducing the memory cost to $O(n + \nnz(A))$.
\end{enumerate}

In \cref{fig:linear_system_comparison}, we observe that the direct method runs out of memory around $n + m = 75$k.
The indirect method baseline runs out of memory earlier at $n = 20$k variables.
This is due to overhead from the custom \texttt{torch.autograd.Function} as well as needing to instantiate not just $M$ but also the gradient $\partial \tilde{x} / \partial M$ as dense matrices.
Note that even though $A$ in our instances has only 1\% nonzero elements, the fill-in caused by the term $A^\top A$ can result in $M$ being very dense.
Our matrix-free adjoint method scales linearly with respect to $n$, requiring less than $300$\,MB even at $n = 50$k.
This constitutes a memory reduction of nearly two orders of magnitude, enabling backpropagation at scales where prior approaches run out of memory.


\section{Details on PSD Projection Backward Pass Experiments (Figures 3 and 4)}
\label{appendix:psd_projection_experiments}

We validate our proposed PSD projection backward rule with two complementary experiments, comparing the following four approaches:
\begin{enumerate}
  \item \textbf{Autograd} backpropagates through the projection \cref{eq:psd-proj} using autograd.
  \item \textbf{Denom. Reg.} replaces $F_{ij}$ in \cref{eq:psd-cotangent} with the expression \cref{eq:denom_reg} that avoids the division by zero for repeated eigenvalues. We report results for $\varepsilon \in \{10^{-3}, 10^{-6}\}$.
  \item \textbf{Newton--Schulz} computes the projection by iterating~\cref{eq:newton-schulz}, and the backward pass is computed by unrolling these Newton--Schulz iterations with autograd.
  \item \textbf{DK (Ours)} uses~\cref{eq:psd-cotangent-dk} to backpropagate through the PSD projection.
\end{enumerate}

The first experiment (\cref{fig:psd_projection_comparison_1}) shows a practical scenario where autograd and denominator regularization fail.
We unroll COSMO with a fixed $\rho$ on a single chance-constrained covariance-steering SDP (double-integrator, $N=10$) for $100$ iterations, starting from a random iterate.
We capture the matrix $S_k$ fed to the PSD projection and the true adjoint cotangent at each iteration via a backward hook, and record the Frobenius norm of every rule's vector-Jacobian product, flagging any non-finite result.
As ADMM drives the projected matrix towards the PSD boundary (i.e., towards satisfying the PSD constraint $S_k \succeq 0$), the eigenvalue gap of the projected matrix collapses below machine precision and autograd returns NaN gradients.
Denominator regularization does not fix the issue, even with different values of~$\varepsilon$.
On the other hand, Newton--Schulz and DK remain finite as they work around the degeneracy of repeated eigenvalues described in~\cref{sec:psd_backward}.

In order to validate the numerical properties of the projection gradients, the second experiment (\cref{fig:psd_projection_comparison_2}) compares the four methods on random symmetric matrices with controlled eigenvalue gaps.
The sampling procedure is defined as follows.
Let $n$ be even and let $\delta > 0$ be the desired eigenvalue gap.
We cluster the eigenvalues around two centers $c_1 = -\delta / 2$ and $c_2 = \delta / 2$.
Defining the spread as $\sigma = \delta \cdot 10^{-3}$, the eigenvalues are sampled via:
\begin{equation}
  \lambda_i = \begin{cases}
    -\delta/2 + u_i, & i = 1, \ldots, n / 2, \\
    +\delta/2 + u_i, & i = n / 2 + 1, \ldots, n,
  \end{cases} \quad u_i \sim \mathcal{U}(-\sigma, \sigma).
\end{equation}
The final sampled matrix is given by $S = Q D(\lambda) Q^\top$, where $Q \in \R^{n \times n}$ is random orthogonal.

Fixing $n = 40$, we sweep the eigenvalue gap $\delta$ from $10^{0}$ down to $10^{-10}$ and report the gradient error relative to a finite-difference reference in the left plot of \cref{fig:psd_projection_comparison_2}.
We show the median and IQR over 100 samples.
In this controlled setting, autograd remains finite, even at small eigenvalue gaps.
Moreover, we observe that denominator regularization carries a large bias, even when $\varepsilon \ll \delta$, because the gaps within each cluster are of order $\sigma = 10^{-3}\delta$ and can still fall below $\varepsilon$, so it is not a robust approach.
Newton--Schulz and DK have low error that does not depend on $\delta$.
Importantly, Newton--Schulz can reach lower errors than DK at the cost of computational speed.
This is evidenced by the right subfigure of \cref{fig:psd_projection_comparison_2}, which shows the wall-clock time of the backward pass approaches for varying matrix dimension $n \in \{10, 20, 40, 80, 160, 320\}$ and $\delta = 10^{-5}$.
We plot the median and IQR over 100 samples.
Newton--Schulz is roughly an order of magnitude slower than the other approaches, due to backpropagating through the iterated matrix multiplications of \cref{eq:newton-schulz}.


\section{Training Details}
\label{appendix:training_details}

\begin{table}[h]
  \centering
  \caption{Training and architecture hyperparameters.}
  \label{tab:training_hyperparams}
  \begin{tabular}{l l}
    \toprule
    Setting & Value \\
    \midrule
    \multicolumn{2}{l}{\emph{Optimization}} \\
    \hquad Optimizer                         & Adam \\
    \hquad Learning rate (cosine schedule)   & $10^{-3} \rightarrow 10^{-4}$ \\
    \hquad Super-batch size $B$              & 5 \\
    \hquad Gradient-norm clip                & 1.0 \\
    \hquad Train / validation instances      & 100 / 20 \\
    \addlinespace
    \multicolumn{2}{l}{\emph{Solver}} \\
    \hquad Linear solve                      & CG \& Jacobi preconditioner \\
    \hquad PSD-projection backward           & DK rule \\
    \addlinespace
    \multicolumn{2}{l}{\emph{Hyperparameter policy} ($\pi_\phi$)} \\
    \hquad Hidden layers (ReLU MLP)          & $32 \times 32$ \\
    \hquad Penalty bounds $[\rho_{\min}, \rho_{\max}]$ & $[10^{-3}, 10^{3}]$ \\
    \hquad Initial $(\rho_{\mathrm{eq}}, \rho_{\mathrm{ineq}}, \rho_{\mathrm{soc}}, \rho_{\mathrm{psd}})$ & $(100, 0.1, 0.1, 0.1)$ \\
    \hquad Initial $\alpha$                   & 1.6 \\
    \addlinespace
    \multicolumn{2}{l}{\emph{Warm-start network}} \\
    \hquad Hidden layers (ReLU MLP)          & $64 \times 64$ \\
    \bottomrule
  \end{tabular}
\end{table}

Our architecture consists of two learnable modules that are trained per problem class: a hyperparameter policy $\pi_\phi$ that predicts the ADMM penalties $\rho^k$ and the relaxation parameter $\alpha^k$ at each iteration, and a warm-start network that predicts the initial iterate $(x^0, s^0, y^0)$.
Both are trained end-to-end with Adam through the unrolled COSMO solver to minimize the self-supervised residual loss~\cref{eq:du-objective}.
All runs use double precision and a single GPU (RTX 4090), and unroll COSMO for $K$ iterations.
The reduced linear system~\cref{eq:cosmo_linear_system_reduced} is solved using the preconditioned CG method, and the PSD-projection backward pass uses the DK rule from~\cref{sec:psd_backward}.
We ablated four different choices of preconditioner (none, Jacobi, block-Jacobi, and a truncated Neumann series with $k = 3$ terms), and found that Jacobi dominates all the others, so we use Jacobi preconditioning for all results.
During training, we use a relative CG tolerance of $10^{-8}$.
For inference, we swept absolute tolerances from $10^{-4}$ to $10^{-12}$ and found no change in ordering, so we report results using a $10^{-4}$ absolute tolerance, which is the fastest for all methods.
Problems are trained in \emph{super-batches} by stacking $B$ instances into a single block-diagonal sparse conic program that can be solved with one CG-based forward/backward pass, which is what makes unrolling large-scale SDPs tractable.
We report results for the checkpoint that has the lowest loss~\cref{eq:du-objective} on the held-out validation set.
The non-learned COSMO baselines use the same unrolled solver but replace $\pi_\phi$ with either tuned constant $(\rho, \alpha)$ or the tuned OSQP-style residual-balancing $\rho$ heuristic, and use a zero (default) initialization, so any speedup is attributable to the learned components alone.
\Cref{tab:training_hyperparams} records the hyperparameters that we use to train each model, and \cref{tab:problem_sizes} reports the per-instance problem sizes along with the unroll horizon $K$ and number of training epochs for each problem class.
Each policy was trained in under 10 GPU hours.

\begin{table}[h]
  \centering
  \caption{Per-instance problem sizes and per-class training settings. $n$ and $m$ are the number of variables and constraints of a single instance, and $\mathrm{nnz}(P)$ and $\mathrm{nnz}(A)$ are the number of nonzeros in the cost and constraint matrices. $K$ is the number of unrolled COSMO iterations during training.}
  \label{tab:problem_sizes}
  \small
  \begin{tabular}{l cccc cc}
    \toprule
    Problem & Variables $n$ & Constraints $m$ & $\mathrm{nnz}(P)$ & $\mathrm{nnz}(A)$ & Unroll $K$ & Epochs \\
    \midrule
    Lyapunov LMI        & 820 & 1640 & 0 & 64040 & 400 & 30 \\
    Max Cut             & 1275 & 1325 & 0 & 1325 & 300 & 60 \\
    Lov\'asz $\vartheta$ & 5050 & 7444--7567 & 0 & 7543--7666 & 100 & 50 \\
    Robust KF           & 1204 & 1200 & 300 & 2700 & 100 & 100 \\
    Robust MPC          & 507 & 1425 & 456 & 4125 & 300 & 60 \\
    Double Integrator   & 3904 & 8363 & 68624 & 8052 & 150 & 60 \\
    Unicycle            & 2943 & 4861 & 100203 & 7200 & 150 & 60 \\
    \bottomrule
  \end{tabular}
\end{table}

\section{Policy Parameterization}
\label{appendix:policy_param}

\paragraph{Hyperparameter policy.}
We use a variant of the policy parameterization proposed by~\citet{Oshin2026Deep}.
At each iteration $k$, the policy predicts a per-constraint penalty vector $\rho^k$ and a scalar relaxation parameter $\alpha^k$.
To keep both parameters in a feasible range, the network outputs unconstrained logits $\ell^k$ that are then bounded using a squashed sigmoid function:
\begin{equation*}
  \rho^k = \exp((\log \rho_{\max} - \log \rho_{\min}) \sigma(\ell_\rho^k) + \log \rho_{\min}) \in (\rho_{\min}, \rho_{\max}), \quad \alpha^k = 2 \sigma (\ell_\alpha^k) \in (0, 2),
\end{equation*}
where $\sigma$ is the logistic sigmoid function.

The policy is recurrent in the logit space.
Two small ReLU MLPs predict additive deltas to the previous iterate's logits:
\begin{equation*}
\begin{aligned}
  \ell_\rho^{k+1} &= \ell_\rho^k + \Delta \ell_\rho^k, \quad \Delta \ell_\rho^k = \mathrm{MLP}_\rho(\xi_\rho^k), \\
  \ell_\alpha^{k+1} &= \ell_\alpha^k + \Delta \ell_\alpha^k, \quad \Delta \ell_\alpha^k = \mathrm{MLP}_\alpha(\xi_\alpha^k).
\end{aligned}
\end{equation*}
The MLPs use as features two sets of residuals produced by the solver at iteration $k$.
The KKT residuals are the first set, and they measure global optimality and feasibility of the iterate $(x^k, s^k, y^k)$:
\begin{equation*}
  r_d^k = \frac{P x^k + q + A^\top y^k}{\max \left( \norm{Px^k}_\infty, \norm{q}_\infty, \norm{A^\top y^k}_\infty \right)} \in \R^n,
  \quad
  r_p^k = \frac{A x^k + s^k - b}{\max \left( \norm{Ax^k}_\infty, \norm{s^k}_\infty, \norm{b}_\infty \right)} \in \R^m.
\end{equation*}
The denominator normalizes the residuals following the OSQP-style rule.
These are the exact quantities used to determine convergence of OSQP and COSMO in practice.

The second set consists of the ADMM dual and primal residuals on the slack variables:
\begin{equation*}
  \delta_d^k = s^k - s^{k - 1}, \quad \delta_p^k = s^k - \tilde{s}^k.
\end{equation*}
The first quantity encodes how fast the iterate is changing (a measure of local optimality), while the second quantity captures how far the projected slack $s^k$ lies from the unprojected iterate $\tilde{s}^k$ (a measure of local feasibility).

Overall, these two sets of residuals carry complementary information.
The KKT residuals are path-independent, global signals that directly determine how far the current iterate is from a primal-dual solution, regardless of how the optimizer arrived there.
Meanwhile, the ADMM residuals are inherently local, as they encode the immediate dynamics of the ADMM iteration.
We find that including both signals is key for enabling the policy to maximize both global and local convergence.

In summary, the policy inputs are given by
\begin{equation*}
  \xi_{\rho, i}^k = \begin{bmatrix}
    \ell_{\rho, i}^k \\
    s_i^k \\
    y_i^k \\
    r_{p, i}^k \\
    \log (\norm{r_{d}^k}_\infty + \epsilon) \\
    \delta_{p, i}^k \\
    \delta_{d, i}^k
  \end{bmatrix}, \quad \xi_{\alpha}^k = \begin{bmatrix}
    \ell_{\alpha}^k \\
    \log (\norm{r_{p}^k}_\infty + \epsilon) \\
    \log (\norm{r_{d}^k}_\infty + \epsilon) \\
    \log (\norm{\delta_p^k}_\infty + \epsilon) \\
    \log (\norm{\delta_d^k}_\infty + \epsilon)
  \end{bmatrix},
\end{equation*}
where $\epsilon = 10^{-12}$ is a small constant to avoid the singularity at $\log(0)$.

Importantly, we adopt an initialization that greatly stabilizes the learning process.
The weights and biases of the last layer of each MLP are initialized as all zeros, resulting in an initial prediction of $\Delta \ell^k = 0$ for all $k$.
This means that the policy initially predicts a fixed $\rho^k = \rho^0$ and $\alpha^k = \alpha^0$ for all iterations, and we set $\rho^0$ and $\alpha^0$ using the COSMO defaults.
Moreover, we let the initial $\rho^0$ and $\alpha^0$ be learnable and train them jointly with the policy.
The policy is therefore initialized so that its predictions are identical to using vanilla COSMO with fixed hyperparameters.
Training then allows the policy to gradually improve over this stable baseline.

Finally, for the non-separable second-order and positive semidefinite cones, the scaled ADMM projection in~\cref{eq:cosmo_projection} couples every row of a cone block.
This means that a single shared penalty parameter per block is required for the projection to remain well-defined.
We therefore apply a \emph{rectification} step that averages the predicted per-row $\rho$ within each SOC/PSD block and broadcasts the mean back over the block.
This allows the per-row policy of~\citet{Oshin2026Deep} to be extended to non-separable cones.
The separable cones, namely the zero and nonnegative cones, are left unmodified with their per-row penalties.

\paragraph{Warm-start network.}
Rather than flattening the variable-size and sparse problem data into a fixed-length vector, our warm start is predicted elementwise by two shared MLPs, mirroring the size-agnostic structure of the policy.
One MLP maps each decision variable's feature vector to its $x^0$ component, and a second maps each constraint's feature vector to its $(s^0, y^0)$ components.
The per-variable features consist of $q_i$, the diagonal $P_{ii}$, and the column $2$-norm and nonzero count of $A_{:,i}$.
The per-constraint features consist of $b_i$, the row $2$-norm and nonzero count of $A_{i,:}$, and a one-hot encoding of the cone type.
A pooled global-context vector $(\log(1 + \norm{q}_\infty), \log(1 + \norm{b}_\infty), \log(1 + n), \log(1 + m))$ is appended to every feature vector to capture the instance's scale and size.
As with the policy, the output layers are zero-initialized, so an untrained warm start equals the zero initialization and training learns a warm-start offset from it.

\begin{figure}[h]
  \centering
  \includegraphics[width=\linewidth]{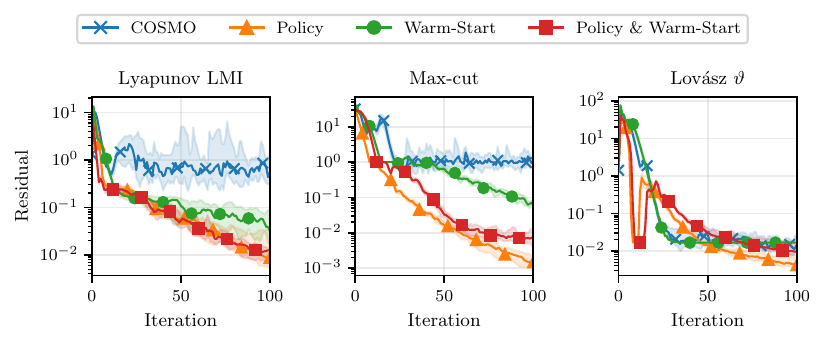}
  \includegraphics[width=\linewidth]{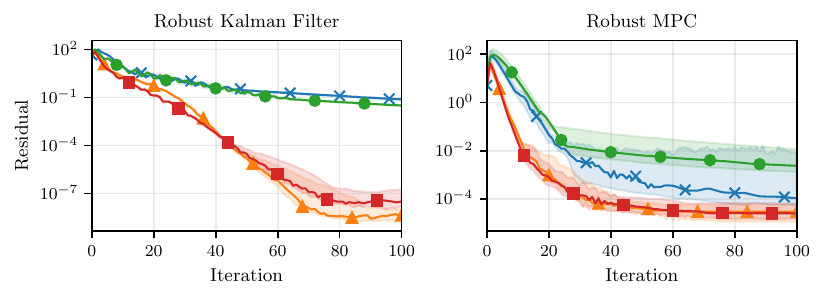}
  \includegraphics[width=\linewidth]{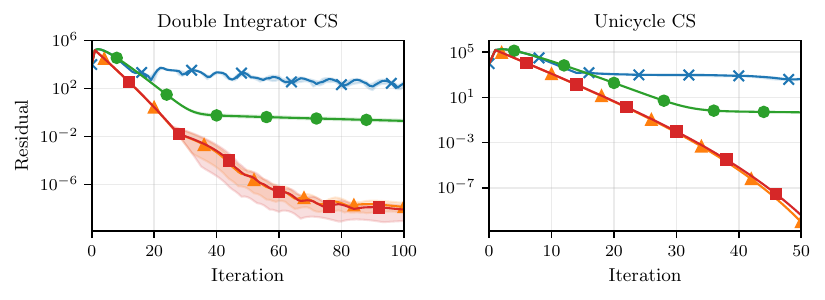}
  \caption{%
    Ablation study comparing the residual convergence of COSMO vs.\ the proposed learned variants. On nearly all of the problem classes, the learned policy exhibits the best performance. The warm start generally delivers an improvement over vanilla COSMO, but occasionally hurts performance, such as on the robust MPC problems. Learning both the policy and warm start jointly does not seem to improve performance over learning the policy alone, so we conclude that most of the performance gain can be attributed to the policy.
  }
  \label{fig:residual_comparison}
\end{figure}

An ablation comparing the performance gains of the warm-start network vs.\ the hyperparameter policy is provided in~\cref{fig:residual_comparison}.
Overall, we find that the hyperparameter policy outperforms the warm start on almost all of the problem classes.
This makes sense because the policy can adapt to local information along the ADMM trajectory, providing a form of closed-loop feedback, whereas the warm start only makes a single open-loop prediction at the start of the optimization.

Lastly, we plot aggregated hyperparameter predictions across all problem classes in~\cref{fig:hyperparameter-predictions-1}.
The learned policy produces qualitatively different behavior compared to the baseline residual-balancing rule.
Moreover, we observe a clear separation between different constraint types that the baseline is unable to reproduce since it only maintains a single scalar $\rho$ (plus a fixed $10^3$ scaling for equality constraints, as in OSQP and COSMO~\citep{Stellato2020OSQP,Garstka2021COSMO}).
The $\alpha$ prediction also differs between problem classes.
Overall, the learned adaptation of the penalty and relaxation parameters at each iteration provides a clear benefit in terms of convergence speed.

\begin{figure}[h]
  \centering
  \includegraphics[width=\linewidth]{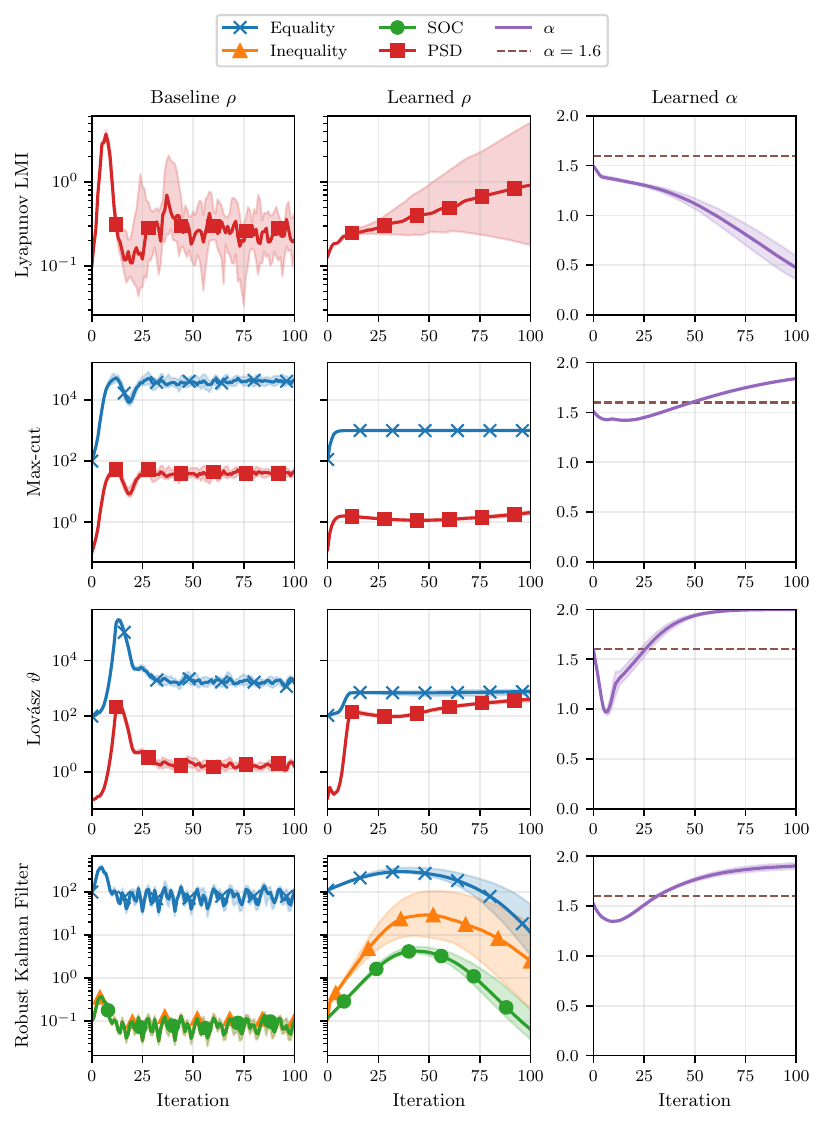}
  \caption{%
    Median and IQR hyperparameter predictions per cone type, with statistics taken over 100 problems. We compare the baseline OSQP residual-balancing rule (left) vs.\ our learned policy (center, right) on all problem classes (rows). The baseline maintains a single scalar $\rho$ that is scaled by $10^3$ for equality constraints~\citep{Stellato2020OSQP}, and uses a fixed $\alpha = 1.6$, denoted by the dashed line on the rightmost plot. The left and center plots share y-axes.
  }
  \label{fig:hyperparameter-predictions-1}
\end{figure}

\begin{figure}[h]
  \ContinuedFloat
  \centering
  \includegraphics[width=\linewidth]{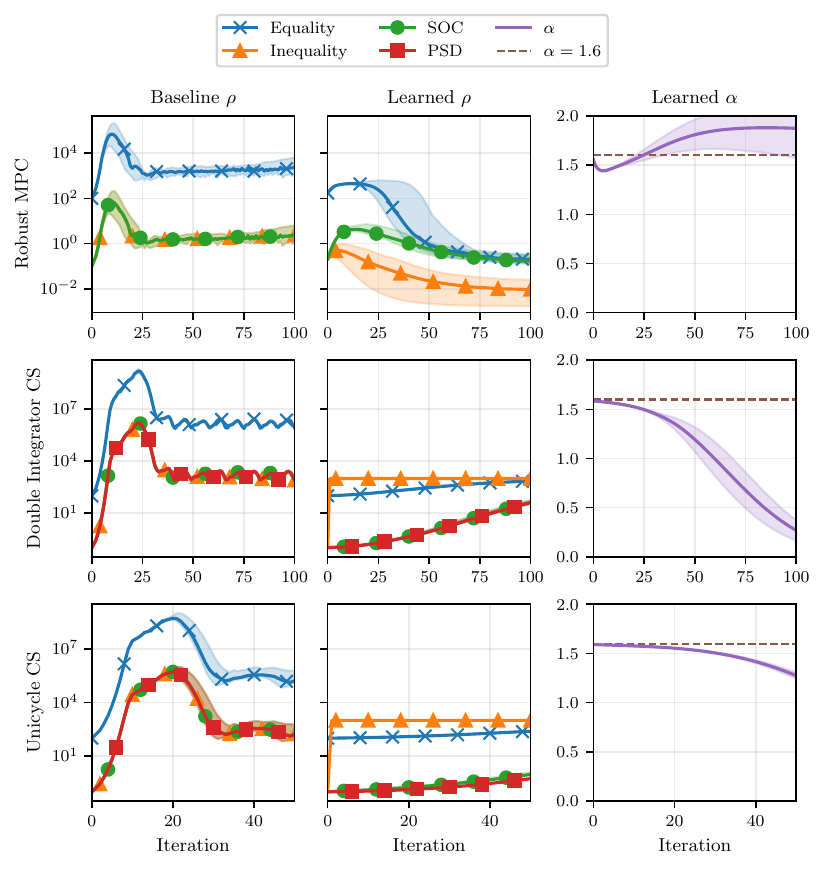}
  \caption{%
    Hyperparameter predictions per cone type (continued).
  }
  \label{fig:hyperparameter-predictions-2}
\end{figure}


\clearpage
\section{Conic Optimization Problems}
\label{appendix:optimization_problems}

This section details the mathematical form of each optimization problem class that we consider in our experiments, and also describes the dataset generation details for each class.
Unless stated otherwise, we sample 100 training and 20 validation instances per class, and evaluate in \cref{tab:cosmo_policies} on a separate set of 100 test instances generated with an independent random seed.

\subsection{Covariance Steering}
Covariance steering seeks a feedback policy that drives the state distribution of a stochastic system from an initial distribution $\mathcal{N}(\mu_0, \Sigma_0)$ to a desired terminal distribution, while minimizing the expected quadratic cost.
We consider the discrete-time stochastic system
\begin{equation*}
  x_{t + 1} = f(x_t, u_t) + w_t, \quad w_t \sim \mathcal{N}(0, W_t), \quad x_0 \sim \mathcal{N}(\mu_0, \Sigma_0),
\end{equation*}
over a horizon of $N$ steps, where $f$ is a (possibly nonlinear) dynamics model and $w_t$ is additive Gaussian process noise.
For a linear time-varying system $x_{t + 1} = A_t x_t + B_t u_t + r_t + w_t$, CS admits a convex reformulation under an affine disturbance-feedback parameterization~\citep{Goulart2006Optimization, Balci2021Covariance}.
Unrolling the recursion from $x_0$ in terms of the state-transition matrix $\Phi(t, \tau) = A_{t - 1} A_{t - 2} \cdots A_\tau$ (with $\Phi(t, t) = I$) gives the closed form
\begin{equation*}
  x_t = \Phi(t, 0)\, x_0 + \sum_{\tau = 0}^{t - 1} \Phi(t, \tau + 1) \left( B_\tau u_\tau + r_\tau + w_\tau \right),
\end{equation*}
which depends only on the inputs and disturbances $u_\tau$ and $w_\tau$ for $\tau < t$.
Stacking the states $X = (x_0, \ldots, x_N)$, the controls $U = (u_0, \ldots, u_{N - 1})$, and the disturbances $W = (w_0, \ldots, w_{N - 1})$, the trajectory rollout can be defined succinctly by the single linear map
\begin{equation*}
  X = \bar{A} x_0 + \bar{B} U + \bar{D} W + \bar{r},
\end{equation*}
where the lifted operators have the block structure
\begin{equation*}
  \bar{A} = \begin{bmatrix} I \\ \Phi(1, 0) \\ \vdots \\ \Phi(N, 0) \end{bmatrix}, \qquad
  \bar{B} = \begin{bmatrix}
    0 & 0 & \cdots & 0 \\
    \Phi(1, 1) B_0 & 0 & \cdots & 0 \\
    \Phi(2, 1) B_0 & \Phi(2, 2) B_1 & \cdots & 0 \\
    \vdots & \vdots & \ddots & \vdots \\
    \Phi(N, 1) B_0 & \Phi(N, 2) B_1 & \cdots & \Phi(N, N) B_{N - 1}
  \end{bmatrix}.
\end{equation*}
$\bar{D}$ shares the same block lower-triangular structure as $\bar{B}$, with each $B_i$ replaced by the disturbance-input gain $D_i$, and the offset $\bar{r}$ stacks the propagated terms $\bar{r}_t = \sum_{\tau < t} \Phi(t, \tau + 1) r_\tau$.

Let $\bar{A}_t, \bar{B}_t, \bar{D}_t$ be the $t$-th block rows of the lifted operators.
We parameterize the control sequence using the affine disturbance feedback form:
\begin{equation*}
  U = K W + V,
\end{equation*}
where $V$ is a feedforward sequence and $K$ is a strictly causal, block lower-triangular feedback gain such that $u_t$ depends only on the past disturbances $w_0, \ldots, w_{t - 1}$.
Substituting this parameterization into the lifted dynamics decouples the per-step mean and covariance, which are given by
\begin{equation*}
  \mu_t = \bar{A}_t \mu_0 + \bar{B}_t V + \bar{r}_t, \qquad
  \Sigma_t = \bar{A}_t \Sigma_0 \bar{A}_t^\top + (\bar{D}_t + \bar{B}_t K) \Sigma_W (\bar{D}_t + \bar{B}_t K)^\top,
\end{equation*}
where $\Sigma_W = \blkdiag(W_0, \ldots, W_{N - 1})$.
The mean is affine in $(\mu_0, V)$ and each covariance is a convex quadratic in $K$, so the expected quadratic cost is convex, and any covariance upper bound $\Sigma_t \preceq M$ becomes a linear matrix inequality via the Schur complement:
\begin{equation*}
  \begin{bmatrix}
    M - \bar{A}_t \Sigma_0 \bar{A}_t^\top & (\bar{D}_t + \bar{B}_t K) \Sigma_W^{1/2} \\
    \Sigma_W^{1/2 \top} (\bar{D}_t + \bar{B}_t K)^\top & I
  \end{bmatrix} \succeq 0.
\end{equation*}
We could instead enforce the PSD constraints with a low-rank factorization $S = L L^\top$ of the PSD matrix $S$, which removes the need for a PSD projection, but we deliberately keep the convex conic form.
The factorization makes the problem nonconvex, so the linear system in \cref{eq:cosmo_linear_system} would become a nonlinear least-squares solve, fundamentally changing the algorithm, and the self-supervised loss based on the KKT residuals would no longer be valid.

For nonlinear $f$, we embed the convex CS problem in a sequential convex programming (SCP) loop inspired by iterative covariance steering~\citep{Ridderhof2019Nonlinear}.
At each outer iteration, we linearize $f$ around the current mean and control reference $(\mu^\text{ref}, u^\text{ref})$ to obtain $(A_t, B_t, r_t)$, solve the resulting convex subproblem, and update the reference under a trust-region ratio test.
A virtual control $\nu_t$ is added to the linearized mean dynamics with an exact $\ell_1$ penalty $\lambda_\nu \norm{\nu_t}_1$, which keeps every subproblem feasible.
Each SCP subproblem is therefore given by the following conic program:
\begin{equation*}
\begin{aligned}
  \minimize_{\mu, V, K, \nu} \hquad & \sum_{t = 0}^{N - 1} \left( \E[x_t^\top Q x_t] + \E[u_t^\top R u_t] + \lambda_\nu \norm{\nu_t}_1 \right) + \E[x_N^\top Q_N x_N], \\
  \st \hquad & \mu_{t + 1} = A_t \mu_t + B_t v_t + r_t + \nu_t, \\
  & \mu_0 = \bar{\mu}_0, \quad \mu_N = \mu_\text{goal}, \\
  & \Sigma_t \preceq \Sigma_\text{max}, \quad \Sigma_N \preceq \Sigma_N^\text{des}, \\
  & \norm{\mu_t - \mu_t^\text{ref}}_2 \leq \rho_x, \quad \norm{v_t - v_t^\text{ref}}_2 \leq \rho_u,
\end{aligned}
\end{equation*}
where the expectations are defined as $\E[x_t^\top Q x_t] = \mu_t^\top Q \mu_t + \tr(Q \Sigma_t)$ (and analogously for the input and terminal costs), and $\rho_x, \rho_u > 0$ are the trust-region radii determined by the outer loop.
After each subproblem solve, we accept or reject the step based on the ratio of actual to predicted cost decrease, expanding or shrinking the trust region accordingly, and declare convergence once the nonlinear residuals (nonlinear dynamics defect, virtual-control norm, etc.) fall below $10^{-3}$.

We consider two systems in our experiments.
The 2D double integrator has state $x = (p_x, p_y, v_x, v_y)$ and control $u = (a_x, a_y)$, with the exact zero-order-hold discretization
\begin{equation*}
  A = \begin{bmatrix} I_2 & \Delta t\, I_2 \\ 0 & I_2 \end{bmatrix}, \quad
  B = \begin{bmatrix} \tfrac{1}{2} \Delta t^2 I_2 \\ \Delta t\, I_2 \end{bmatrix}, \quad
  W = \diag\!\left( \tfrac{\Delta t^3}{3} \sigma_a^2 I_2, \ \Delta t\, \sigma_a^2 I_2 \right),
\end{equation*}
with acceleration-noise scale $\sigma_a = 0.1$.
We use a time horizon of $N = 30$ with a discretization step size of $\Delta t = 0.5$.
The dynamics are linear, so the SCP loop converges in a single iteration, and each instance is a single convex subproblem.
The initial mean for each instance is fixed at $\bar{\mu}_0 = 0$ and the target mean is sampled from $\mathcal{U}(-\bar{x}, \bar{x})$, where $\bar{x} = (5, 5)$.
The starting covariance is given as $\Sigma_0 = \sigma_0^2 I$, with $\sigma_0 = 0.01$ for the training and validation instances and $\sigma_0 \sim \mathcal{U}(0.005, 0.02)$ for the test instances, and the terminal covariance is given as $\Sigma_N^\text{des} = \sigma_N^2 I$ with $\sigma_N = 0.25$.
The cost function parameters are given by $Q = \diag(1, 1, 0.3, 0.3)$, $R = 0.1 I_2$, and $Q_N = 10 Q$.

The kinematic unicycle has state $x = (p_x, p_y, \theta)$ (position and heading) and control $u = (v, \omega)$ (linear and angular velocity), with continuous-time dynamics
\begin{equation*}
  \dot{x} = f(x, u) = \begin{bmatrix}
    v \cos\theta \\
    v \sin\theta \\
    \omega
  \end{bmatrix}.
\end{equation*}
We discretize the dynamics using the explicit Euler step $x_{t + 1} = x_t + f(x_t, u_t) \Delta t$.
The process-noise covariance is given by $W = \diag(\sigma_p^2, \sigma_p^2, \sigma_\theta^2) \Delta t$ with $(\sigma_p, \sigma_\theta) = (0.02, 0.01)$.
The nonlinear coupling from the heading $\theta$ makes each linearization point distinct, so the SCP loop requires several outer iterations per instance, around 10--15 in our testing.
We fix the starting mean at the origin and sample its heading uniformly from $[-\pi, \pi]$.
The target mean is sampled from $\mathcal{U}(-\bar{x}, \bar{x})$, where $\bar{x} = (5, 5, \pi)$.
The boundary standard deviations are sampled using $\sigma_0 \sim \mathcal{U}(0.01, 0.05)$ and $\sigma_N \sim \mathcal{U}(0.05, 0.25)$.
We use a time horizon of $N = 30$ with a discretization step size of $\Delta t = 0.1$, along with cost function parameters $Q = \diag(1, 1, 0.1)$, $R = \diag(0.1, 0.05)$, $Q_N = 10 Q$, and $\lambda_\nu = 10^4$.
The trust region is initialized at $\rho_x = \rho_u = 5$ and adapted multiplicatively (shrinking by $0.5$, expanding by $2.0$) each outer iteration.
We sample 12 nonlinear CS problems and run 10 SCP iterations on each, solving every subproblem with Clarabel, to generate 120 subproblem instances.
The 100 subproblems from the first 10 scenarios are used for training and the 20 from the remaining 2 for validation.
The 100 test subproblems used in \cref{tab:cosmo_policies} come from 10 further scenarios sampled with an independent seed, so no subproblem from a test scenario appears in the training set.
We sample 20 new scenarios for the evaluation in \cref{fig:nonlinear_cs_comparison}.

\subsection{Lyapunov Linear Matrix Inequality}
The first class of standalone SDPs we consider is the minimum-trace Lyapunov LMI SDP.
The goal is to compute a Lyapunov certificate for a linear time-invariant system $\dot{x} = Ax$, where $A \in \R^{n \times n}$ is Hurwitz.
Finding $P \succ 0$ such that $\dot{V} = x^\top (A^\top P + P A) x < 0$ for $V(x) = x^\top P x$ shows the origin is globally exponentially stable.
The optimization problem is given as:
\begin{equation*}
\begin{aligned}
  \minimize_{P \in \mathbb{S}^n} \hquad & \tr(P), \\
  \st \hquad & A^\top P + P A \preceq -I, \\
  & P \succeq I.
\end{aligned}
\end{equation*}
The first constraint is the Lyapunov inequality that certifies stability, and we adopt a normalized form such that the problem has a unique minimizer.
The second constraint normalizes the scale of $P$ and guarantees positive definiteness.
Minimizing the trace of $P$ finds the tightest such certificate.
This form of optimization problem appears commonly in optimal control and robotics~\citep{Tedrake2023Underactuated}, as it is fundamental for computing LQR value functions, invariant sets, etc.

We generate instances with $n = 40$.
To avoid overfitting to a particular spectrum, each $A$ is drawn with equal probability from one of two families.
The first is $A = Q T Q^\top$, where $Q$ is random orthogonal and $T$ is upper-triangular with diagonal elements $T_{ii} \sim \mathcal{U}(-2, -0.1)$ and off-diagonal elements $T_{ij} \sim \mathcal{N}(0, 0.1^2)$.
The small off-diagonal scale limits the non-normality of $A$, which would otherwise make the Lyapunov operator severely ill-conditioned.
The second is $A = Q R Q^\top$, where $Q$ is random orthogonal and $R = \blkdiag(R_1, R_2, \ldots, R_{n/2})$ is a block-diagonal matrix of damped oscillators:
\begin{equation*}
  R_i = \begin{bmatrix}
    -\zeta_i \omega_i & \omega_i \\
    -\omega_i & -\zeta_i \omega_i
  \end{bmatrix},
\end{equation*}
with natural frequency $\omega_i \sim \mathcal{U}(0.5, 5)$ and damping ratio $\zeta_i \sim \mathcal{U}(0.05, 0.5)$.
We resample any $A$ whose spectral abscissa is not below $-0.05$.

\subsection{Max Cut}
The max-cut SDP is a fundamental problem in graph theory and optimization~\citep{Goemans1995Improved}.
For an undirected weighted graph $G$, the Laplacian is defined as $L = D - A$, where $D$ is the degree matrix and $A$ is the adjacency matrix.
The max-cut SDP is written as
\begin{equation*}
  \begin{aligned}
  \maximize_{X \in \mathbb{S}^n} \hquad & \frac{1}{4} \langle L, X \rangle, \\
  \st \hquad & X_{ii} = 1, \quad \forall i = 1, \ldots, n, \\
  & X \succeq 0.
  \end{aligned}
\end{equation*}
We sample random Erd\H{o}s--R\'enyi graphs with $n = 50$ nodes and edge probability $p = 0.5$.

\subsection{Lov\'asz $\vartheta$}

The Lov\'asz $\vartheta$ SDP is a classical problem in graph theory that seeks an upper bound on the Shannon capacity of a graph~\citep{Lovasz1979Shannon}.
Given a graph $G = (V, E)$, the problem is given by
\begin{equation*}
  \begin{aligned}
  \maximize_{X \in \mathbb{S}^n} \hquad & \langle J, X \rangle, \\
  \st \hquad & \tr(X) = 1, \\
  & X_{ij} = 0, \quad \forall (i, j) \in E, \\
  & X \succeq 0,
  \end{aligned}
\end{equation*}
where $J$ is a matrix of all ones. Similar to the max-cut problems, we sample random Erd\H{o}s--R\'enyi graphs with $n = 100$ nodes and edge probability $p = 0.5$.

\subsection{Robust Kalman Filter}
Robust Kalman filtering is a variant of the traditional Kalman filter that is more robust to large outliers in the sensor measurements.
This second-order cone program is given by
\begin{equation*}
  \begin{aligned}
    \minimize \hquad & \sum_{t = 0}^{N - 1} \left( \norm{w_t}_2^2 + \tau \phi_\delta(v_t) \right), \\
    \st \hquad & x_{t + 1} = Ax_t + Bw_t, \\
    & y_t = Cx_t + v_t,
  \end{aligned}
\end{equation*}
where $\phi_\delta$ with $\delta > 0$ is the Huber function:
\begin{equation*}
  \phi_\delta(v) = \begin{cases}
    \norm{v}_2^2 & \text{if}\ \norm{v}_2 \leq \delta, \\
    \delta (2 \norm{v}_2 - \delta) & \text{if}\ \norm{v}_2 > \delta.
  \end{cases}
\end{equation*}

Our instances are adapted from the CVXPY robust Kalman filtering example~\citep{cvxpy_robust_kalman}. We use a double integrator in 2D with linear drag:
\begin{equation*}
  A = \begin{bmatrix}
    1 & 0 & (1 - \tfrac{1}{2} \gamma \Delta t) \Delta t & 0 \\
    0 & 1 & 0 & (1 - \tfrac{1}{2} \gamma \Delta t) \Delta t \\
    0 & 0 & 1 - \gamma \Delta t & 0 \\
    0 & 0 & 0 & 1 - \gamma \Delta t
  \end{bmatrix}, \quad B = \begin{bmatrix}
    \tfrac{1}{2} \Delta t^2 & 0 \\
    0 & \tfrac{1}{2} \Delta t^2 \\
    \Delta t & 0 \\
    0 & \Delta t
  \end{bmatrix}, \quad C = \begin{bmatrix}
    1 & 0 & 0 & 0 \\
    0 & 1 & 0 & 0
  \end{bmatrix},
\end{equation*}
with $\Delta t = 0.1$ and drag coefficient $\gamma = 0.05$, and we set $\tau = \delta = 2$.
We sample process noise $w_t \sim \mathcal{N}(0, 1)$ and roll out the dynamics for $N = 100$ steps starting from the origin to generate the state trajectory.
The observations are generated using the measurement model $y_t = Cx_t + v_t$, with $v_t$ sampled via
\begin{equation*}
  v_t \sim \begin{cases}
    \mathcal{N}(0, I) & \text{with probability}\ p = 0.8, \\
    \mathcal{N}(0, \sigma_\text{outlier}^2 I) & \text{otherwise},
  \end{cases}
\end{equation*}
where $\sigma_\text{outlier} = 20$.

\subsection{Robust Model Predictive Control}

Finally, we consider robust MPC problems for discrete-time linear systems of the form $x_{t + 1} = A x_t + B u_t + w_t$, where $w_t \in \mathcal{W}$ is a bounded disturbance in the ellipsoid $\mathcal{W} = \{w \mid w^\top \Sigma_w^{-1} w \leq 1 \}$ for $\Sigma_w \succ 0$.
Let $A_x x_t \leq b_x$ and $A_u u_t \leq b_u$ be the affine state and control constraints.
The goal in robust MPC is to determine a nominal trajectory that satisfies the constraints under all possible realizations of uncertainty.
The covariance of each state can be computed offline through $\Sigma_{t + 1} = A \Sigma_t A^\top + \Sigma_w$ and $\Sigma_0 = 0$.
The half-space constraints of the form $a_i^\top x \leq b_i$ are then tightened using $a_i^\top x \leq b_i - \kappa \lVert\Sigma_t^{1/2} a_i\rVert_2$, where $\kappa > 0$ is the confidence parameter.
Letting $[d_t]_i = \norm{\Sigma_t^{1/2} a_i}_2$, the robust MPC problem is given as
\begin{equation*}
\begin{aligned}
  \minimize \hquad & \sum_{t = 0}^{N - 1} (x_t^\top Q x_t + u_t^\top R u_t + \gamma \tau_t) + x_N^\top Q_N x_N + \gamma_\text{term} \tau_\text{term}, \\
  \st \hquad & x_{t + 1} = A x_t + B u_t, \\
  & x_0 = x_\text{init}, \\
  & A_x x_t \leq b_x - \kappa d_t, \\
  & A_u u_t \leq b_u, \\
  & \norm{u_t}_2 \leq \tau_t, \\
  & \norm{x_N - x_\text{target}}_2 \leq \tau_\text{term},
\end{aligned}
\end{equation*}
where $\tau_t$ and $\tau_\text{term}$ are the epigraph variables for the control effort and terminal constraint, respectively, with $\gamma, \gamma_\text{term} > 0$.

We generate problem instances with $N = 50$ for a chain of three unit masses connected in series by springs and dampers, with the first mass anchored to a wall and the last mass free.
Each mass is actuated by an independent force, so the state $x = (p, v) \in \R^6$ stacks the mass positions and velocities and the control is $u \in \R^3$.
The continuous-time dynamics are
\begin{equation*}
  \dot{x} = \begin{bmatrix}
    0 & I_3 \\
    -K & -c I_3
  \end{bmatrix} x + \begin{bmatrix}
    0 \\
    I_3
  \end{bmatrix} u, \qquad
  K = k \begin{bmatrix}
    2 & -1 & 0 \\
    -1 & 2 & -1 \\
    0 & -1 & 1
  \end{bmatrix},
\end{equation*}
with spring stiffness $k = 1$ and damping coefficient $c = 0.3$, which we discretize with an exact zero-order hold at $\Delta t = 0.2$ to obtain $(A, B)$.
We use cost matrices $Q = Q_N = I_6$ and $R = I_3$.
The state and control constraints are the symmetric boxes $\lvert p_i \rvert \leq 5$, $\lvert v_i \rvert \leq 1$, and $\lvert u_i \rvert \leq 0.5$.
The disturbance covariance is set to $\Sigma_w = \sigma^2 I$ with $\sigma = 0.01$, and the confidence parameter is $\kappa = 2$.
The dynamics are fixed across instances.
The initial positions and velocities are sampled from $p_\text{init} \sim \mathcal{U}(-1, 1)$ and $v_\text{init} \sim \mathcal{U}(-0.3, 0.3)$ elementwise, and the terminal condition is set as $x_\text{target} = 0$. The penalty parameters are sampled as $\gamma \sim \mathcal{U}(0.01, 0.1)$ and $\gamma_\text{term} \sim \mathcal{U}(1, 10)$.
Some sampled initial conditions admit no trajectory satisfying the tightened constraints, so we solve each candidate with Clarabel and reject it if it is infeasible or if its optimal dual variable is badly scaled ($\norm{y^\star}_2 > 150$); roughly 10\% of candidates are rejected.


\section{Additional Benchmarks}
\label{appendix:additional_benchmarks}

We additionally benchmark our PyTorch implementation against COSMO.jl v0.8.9 (called from Python via \texttt{juliacall}) and Clarabel v0.11.1 on the same instances from \cref{tab:cosmo_policies}.
The wall-clock times to convergence are reported in \cref{tab:additional_benchmarks}.
COSMO.jl and Clarabel run on the CPU with 8 BLAS threads.
For COSMO.jl, $\varepsilon_\mathrm{abs} = 10^{-3}$ and $\varepsilon_\mathrm{rel} = 0$ with termination checked every iteration, as in our solver.
For Clarabel, absolute and relative gap and feasibility tolerances are set to $10^{-3}$.
Chordal decomposition is disabled in both, in order to isolate our contributions.
Every other setting is left at the default (COSMO.jl: Ruiz scaling, adaptive $\rho$, and Anderson acceleration; Clarabel: presolve and equilibration).
Timings for each solver are recorded after an initial warm-up solve.

Overall, COSMO.jl achieves the lowest solve times on the Lyapunov LMI and max-cut SDP classes, while Clarabel is fastest on the second-order cone programs (robust KF and robust MPC).
Our learned solver is fastest on the largest problem classes: Lov\'asz $\vartheta$ and the CS double integrator and unicycle problems (see \cref{tab:problem_sizes} for variable and constraint counts).
Since these classes are substantially larger than the others, they provide evidence that our approach scales effectively to large-scale conic optimization problems.

\begin{table}[h]
  \centering
  \caption{%
    Median [IQR] wall-clock time per instance (in seconds), over 100 held-out test instances.
    Each solver uses a convergence tolerance of $10^{-3}$.
    \emph{Constant}, \emph{Adaptive}, and \emph{Learned} use our PyTorch reimplementation of COSMO with the corresponding policies from \cref{tab:cosmo_policies}.
  }
  \label{tab:additional_benchmarks}
  \tiny
  \begin{tabular}{l ccccc}
    \toprule
    Problem Class & Constant & Adaptive & Learned & COSMO.jl & Clarabel \\
    \midrule
    \hquad Lyapunov LMI & 1.58 [1.04, 3.11] & 2.42 [1.26, 5.39] & 1.09 [0.761, 1.52] & \textbf{0.462 [0.398, 0.590]} & 0.593 [0.519, 0.721] \\
    \hquad Max Cut & 0.396 [0.376, 0.421] & 1.08 [0.788, 2.01] & 0.159 [0.149, 0.179] & \textbf{0.0411 [0.0382, 0.0492]} & 0.229 [0.205, 0.239] \\
    \hquad Lov\'asz $\vartheta$ & -- & 1.01 [0.547, 12.2] & \textbf{0.359 [0.333, 0.391]} & 0.599 [0.523, 0.773] & 4.27 [4.23, 4.31] \\
    \hquad Robust KF & 4.62 [4.50, 4.75] & 4.64 [4.52, 4.79] & 0.884 [0.831, 0.966] & 0.0169 [0.0167, 0.0171] & \textbf{0.00208 [0.00205, 0.00213]} \\
    \hquad Robust MPC & 1.02 [0.842, 1.35] & 0.821 [0.633, 1.02] & 0.216 [0.196, 0.248] & 0.0137 [0.0134, 0.0154] & \textbf{0.00193 [0.00173, 0.00195]} \\
    \hquad Double Integrator & 0.971 [0.531, 1.57] & 4.81 [3.03, 9.20] & \textbf{0.633 [0.428, 0.734]} & 2.73 [2.53, 2.95] & 14.1 [12.8, 14.5] \\
    \hquad Unicycle & 31.4 [21.3, 40.9] & -- & \textbf{0.911 [0.808, 0.992]} & 1.48 [1.22, 2.16] & 3.82 [3.53, 3.98] \\
    \bottomrule
  \end{tabular}
\end{table}


\section{Control Performance and Failure Modes}
\label{appendix:control_performance}

The iteration counts and wall-clock times in the main text measure solver performance, so this section validates that the solutions returned by the learned solver also produce good control.

\paragraph{Failure modes for nonlinear covariance steering.}
In the \cref{fig:nonlinear_cs_comparison} experiment, learned COSMO fails on 4 of the 20 scenarios and Clarabel fails on 8 of the 20.
Clarabel also fails on every scenario that learned COSMO fails on.
Each inner COSMO solve is capped at 1000 iterations and 20\,s, but the iteration cap is never reached in practice: within the 20\,s limit, baseline COSMO runs roughly 450 ADMM iterations.
A run fails in one of two ways: either the inner subproblems are not solved, or they are solved but the outer SCP loop stalls with the max SCP residual above the $10^{-3}$ tolerance.
For learned COSMO, on one scenario, every subproblem reaches the 20\,s time limit of the inner solver; the SCP residual still decreases, but is only $3.3 \times 10^{-3}$ by the time the budget of 20 SCP iterations runs out.
A second scenario ends at the same residual after a steady decrease, with all but one subproblem converged.
In the remaining two scenarios, the first subproblem reaches the 20\,s limit, and the max SCP residual stalls at $0.29$ and $0.99$ despite every later subproblem converging in under 80 iterations.
For Clarabel, on five scenarios, at every SCP iteration the interior-point method stops after 6--10 iterations without returning a solution, so the reference trajectory is never updated.
On the other three scenarios, every subproblem is solved, but the SCP residual stalls between $7.6 \times 10^{-3}$ and $6.7 \times 10^{-2}$.
Inexact subproblem solutions do not prevent convergence: learned COSMO solves three scenarios where subproblems hit the time limit.

\paragraph{Monte Carlo validation of the covariance bounds.}
We validate the policies returned by SCP on the true nonlinear system.
For each solved scenario, we roll out the affine disturbance-feedback policy $u_t = v_t + \sum_{\tau < t} K_{t, \tau} w_\tau$ on the Euler-discretized unicycle under 2000 sampled realizations of the initial state $x_0 \sim \mathcal{N}(\mu_0, \Sigma_0)$ and process noise $w_t \sim \mathcal{N}(0, W)$, with the same samples for both solvers.
On every solved scenario, for both solvers, the empirical terminal covariance $\hat{\Sigma}_N$ satisfies $\hat{\Sigma}_N \preceq \Sigma_N^\text{des}$ with margin: $\lambda_\text{max}(\hat{\Sigma}_N - \Sigma_N^\text{des})$ lies between $-5.0 \times 10^{-2}$ and $-5.1 \times 10^{-4}$, with median $-1.2 \times 10^{-2}$, in agreement with the value predicted by the linearized model.
The empirical terminal mean is within $7.1 \times 10^{-3}$ of the goal, and the median ratio of the empirical expected cost to the cost of the SCP solution is $1.0006$ for learned COSMO and $1.0003$ for Clarabel.

\paragraph{Closed-loop robust MPC.}
We run the robust MPC controller from \cref{appendix:optimization_problems} on 10 closed-loop rollouts of 50 steps each.
Each episode samples $x_\text{init}$, $\gamma$, and $\gamma_\text{term}$ from the same distribution as the test dataset.
At each step, the control problem (horizon $N = 50$) is re-solved from the current state and the first control is applied to the linear dynamics.
Process noise $w_t \sim \mathcal{N}(0, \Sigma_w)$ is added, with the same noise sequence applied for each solver.
We compare learned COSMO, the adaptive COSMO baseline, and Clarabel, each with a convergence tolerance of $10^{-3}$.
Both COSMO arms have a time budget of 60\,s per MPC step.
When a solver does not converge, its returned iterate is applied (zero control if the iterate is non-finite or no solution is returned).

On 9 of the 10 episodes, all three solvers converge at every step and no state constraint is violated.
The closed-loop cost of learned COSMO is within $0.1\%$ of that of Clarabel on all but one episode, where it is $0.9\%$ lower.
Learned COSMO requires a median of 12 iterations ($0.29$\,s) per step, compared with 50 iterations ($2.10$\,s) for COSMO, while Clarabel on the CPU takes $0.035$\,s per step at this problem size.
On the remaining episode, all three solvers fail within the first few steps (COSMO at step 5, learned COSMO and Clarabel at step 7).
Learned COSMO reaches the 60\,s limit with a non-finite iterate and Clarabel returns no solution, so both apply zero control and follow the same trajectory, which violates the state box constraints on four steps by up to $0.42$.
COSMO also reaches the 60\,s limit, but applies its inexact iterates, which violate the state constraints on three steps by up to $0.08$.


\section{Details on Batched Solve Experiment (Table 2)}
\label{appendix:batched_experiment}

\paragraph{Instances.}
We sample new batches of $B \in \{8, 16, 32, 64\}$ instances from the same generator and parameters as the held-out test sets from~\cref{tab:cosmo_policies} (max cut: $n = 50$ nodes with edge probability $0.5$; Lov\'asz $\vartheta$: $n = 100$ nodes with edge probability $0.5$).
Each batch is seeded by its batch size and sample index, so no batch repeats a training, test, or other batch instance, and the same batch is given to all three solvers.
\Cref{tab:batched_solves} reports $B = 16$ and $B = 64$; over $B = 8, 16, 32, 64$, the ratio of median times, COSMO.jl over learned, is $1.4\times$, $2.9\times$, $5.0\times$, $8.1\times$ on max cut and $5.6\times$, $6.3\times$, $9.4\times$, $14\times$ on Lov\'asz $\vartheta$.

\paragraph{Learned solver.}
The learned solver uses the hyperparameter policy and warm-start network of~\cref{tab:cosmo_policies}, with $\rho$ updated at every iteration as during training.
The $B$ instances are stacked into one block-diagonal conic program, in the same fashion as the training super-batches of~\cref{appendix:training_details}.
Every iteration acts on the whole super-batch: one CG solve, one projection per cone, and one policy evaluation, which is row-wise and therefore assigns each instance its own penalties.
The only per-instance operation is the convergence check, which reduces the primal and dual residuals over each instance's rows at every iteration.
The solve terminates when every instance satisfies $\varepsilon = 10^{-3}$.
Instances that have already converged keep iterating, so the batch's iteration count is that of its slowest instance (the median count rises from 110 to 122 on max cut and from 154 to 434 on Lov\'asz $\vartheta$ between $B = 8$ and $B = 64$).
The linear system is solved with Jacobi-preconditioned CG at an absolute tolerance of $10^{-4}$, and each batch solve is capped at 10000 iterations and 60\,s.
The learned arm solved all 10 batches at every $B$ on both classes without reaching the budget.

\paragraph{Timing.}
All solvers run on the same machine (i9-13900K CPU, RTX 4090 GPU, double precision), and each is run once untimed on a warm-up batch, which absorbs the CUDA initialization and Julia's JIT compilation time.
The learned arm's time is the CUDA-synchronized wall-clock time from the start of super-batch assembly to the end of the solve (assembly and solver construction account for 1--3\% of the total time).
Each baseline's time is the sum over the batch of the solver's reported total time (COSMO.jl's \texttt{solver\_time}, which includes setup; Clarabel's \texttt{solve\_time}, which includes equilibration and factorization), i.e., the total cost of a native sequential loop.


\section{Details on Out-of-Distribution Experiment (Table 3)}
\label{appendix:ood_experiment}

\paragraph{Instances.}
Each class is scaled using a size parameter (number of nodes $n$ for max cut, planning horizon $N$ for robust KF and robust MPC), with every other parameter kept at its training configuration (\cref{appendix:optimization_problems}).
These classes were chosen because a single parameter controls the problem size: the number of variables grows quadratically with $n$ for max cut and linearly with $N$ for robust KF and robust MPC.
Scaling $n$ enlarges the PSD block of the max-cut SDP from $50 \times 50$ to $500 \times 500$ (from 1275 to 125250 variables).
Scaling the horizon for robust KF and robust MPC appends extra second-order cones of the same size, so the problems grow in the number of blocks rather than in the block size.
We generate 20 held-out test instances per size.

\paragraph{Arms.}
\emph{Learned} is the class's hyperparameter policy and warm-start network, applied without retraining or any architectural change: the policy acts row-wise on per-constraint residuals and iterates (\cref{appendix:policy_param}) and the warm-start network likewise maps per-row features, so neither has parameters tied to $n$, $m$, or the number of cones, and the same weights apply at every size.
On the SOCPs, the predicted $\rho$ is adopted at every iteration.
On max cut, the predicted $\rho$ is adopted only when it differs from the current value by more than a factor of 5, which is the same gate that the residual-balancing baseline uses.
The predicted penalties diverged when adopted at every iteration, which we attribute to residual scales shifting with problem size, so overly frequent feedback destabilizes the policy.
\emph{Baseline} is COSMO with the OSQP-style residual-balancing $\rho$ rule at the same adoption gate of 5 and the default zero initialization, i.e., the \emph{Adaptive} arm of~\cref{tab:cosmo_policies}.
The \emph{Constant} arm is omitted because its tuned $\rho$ is specific to the training size.
Both arms run in our PyTorch solver on the GPU with the same Jacobi-preconditioned CG solve (absolute tolerance $10^{-4}$), so the comparison isolates the learned components.

\paragraph{Protocol.}
After an untimed warm-up solve, each arm performs one timed solve per instance to $\varepsilon = 10^{-3}$ tolerance on the ADMM residuals, with convergence checked every iteration.
The reported time is the CUDA-synchronized wall-clock time of the solve.
A solve that does not converge within 10000 iterations or 600\,s is charged the budget caps, following the OSQP benchmark protocol~\citep{Stellato2020OSQP}, but no solve of either arm failed at any size.
\Cref{tab:ood_generalization} reports the median [IQR] of iterations and time per arm, and the speedup as the median over instances of the per-instance baseline-to-learned ratio, which is why it can differ from the ratio of the two medians.


\end{document}